\documentclass[review]{elsarticle}
\RequirePackage[top=1.5cm,bottom=1.5cm,left=2.cm,right=2.cm,includehead,includefoot]{geometry}
\usepackage{amssymb}
\usepackage{fancyhdr}
\usepackage{amsmath,amsfonts,amsthm}
\usepackage{mathrsfs}
\usepackage{amsfonts}
\usepackage{cases}
\usepackage{verbatim}
\usepackage{graphicx}      % needed for including graphics
\usepackage{epstopdf}
\usepackage{grffile}
\graphicspath{{figures/}}
\usepackage{subfig}
\usepackage{caption}
\usepackage{multirow}
\usepackage{threeparttable}
\usepackage{bm}
\usepackage{makecell}
\usepackage{booktabs}
\usepackage{hyperref}
\usepackage{cleveref}

\newtheorem{thm}{Theorem}[section]

\newtheorem{lemma}[thm]{Lemma}

\theoremstyle{remark}

\usepackage[frak=mma]{mathalfa}

\newcommand{\la} \langle
\newcommand{\ra} \rangle
\begin{document}

\begin{frontmatter}



\title{Render unto Numerics: Orthogonal Polynomial Neural Operator for PDEs with Nonperiodic Boundary Conditions\tnoteref{mytitlenote}}
%\tnotetext[mytitlenote]{This work was supported by the National Natural Science Foundation of China (Grant No. 11571366), the National Natural Science Foundation of Hunan (Grant No. S2017JJQNJJ0764), Research Fund of NUDT (Grant No. ZK17-03-27), and the fund from Hunan Provincial Key Laboratory of Mathematical Modeling and Analysis in Engineering (dGrant No. 2018MMAEZD004).}


\author[nudt]{Liu Ziyuan}
\author[nudt]{Wang Haifeng}
\author[nudt]{Bao Kaijun}
\author[nudt]{Qian Xu *}
\ead{qianxu@nudt.edu.cn}
\author[nudt]{Zhang Hong}
\author[nudt,hpc]{Song Songhe}

\address[nudt]{College of Science, National University of Defense Technology, Changsha 410073, China}
\address[hpc]{State Key Laboratory of High Performance Computing, National University of Defense Technology, Changsha 410073, China}

\begin{abstract}
  By learning the mappings between infinite function spaces using carefully designed neural networks, the operator learning methodology has exhibited significantly more efficiency than traditional
methods in solving complex problems such as differential equations, but faces concerns about their accuracy and reliability. To overcomes these limitations, combined with the structures of the spectral numerical method, a general
neural architecture named spectral operator learning (SOL) is introduced, and one variant called the orthogonal polynomial neural operator (OPNO), developed for PDEs with Dirichlet,
Neumann and Robin
boundary conditions (BCs), is proposed later. The strict BC satisfaction properties and the universal approximation capacity of the OPNO are theoretically proven. A variety of numerical experiments with physical backgrounds show that the OPNO outperforms other existing deep learning methodologies, as well as the traditional 2nd-order finite difference
method (FDM) with a considerably fine mesh (with the relative errors reaching the order of $10^{-6}$), and is up
to almost $5$ magnitudes faster than the traditional method.

\end{abstract}

\begin{keyword}
deep learning-based PDE solver, neural operator, spectral method, AI4science, scientific machine learning.

\MSC 47-08 \sep 65D15 \sep 65M22 \sep 68Q32 \sep 68T07

\end{keyword}

\end{frontmatter}


\section{Introduction}\label{Introduction}
Differential equations are the foundational models in numerous fields of modern science and engineering, and for decades, their solving processes have been dominated by numerical methods. However, recent research
has shown that deep neural networks have an extraordinary capacity to solve highly nonlinear problems
and the potential to develop algorithms more efficient than numerical methods. On the other hand, the stability and accuracy of numerical methods are ensured by numerical approximation theory, and these methods often lead to sparse systems, all of which are extremely desirable
in the deployment of deep neural networks. Therefore, the following question arises: how can we incorporate these features to propose an innovative deep learning-based framework for partial differential equation (PDE) solvers?

In recent years, significant efforts have been oriented towards the development of deep learning-based PDE solvers. Some approaches focus on directly approximating PDE solutions with neural networks, such as the deep Galerkin method \cite{sirignano2018dgm}, the deep Ritz
method \cite{yu2018deep,liao2019deep} and physically informed neural networks \cite{raissi2019physics,wang2022respecting} and are able to overcome the curse of dimensionality in theory but need retraining
every time the PDE parameters or conditions are slightly changed. Furthermore, operator
learning, a general methodology that learns the solution mappings between input and output function spaces, can efficiently generate the solution of an entire
family of PDEs and is \textit{roughly} categorized into two families of methods: FNO-like operators such as Fourier neural
operators itself \cite{li2021fourier}, the multiwavelet-based
neural operator \cite{gupta2021multiwavelet}, the Galerkin transformer \cite{cao2021choose}, spectral neural operators \cite{fanaskov2022spectral}, and the integral autoencoder \cite{ong2022iae}; or DeepONet-like
operators such as DeepONet \cite{lu2019deeponet}, DeepM\&Mnet \cite{cai2021deepm}, and MIONet \cite{jin2022mionet}. Among the available operator learning methods, FNOs are
both efficient and accurate and therefore may be the most fruitful and promising neural operator method to date for practical engineering applications, such as high-resolution weather forecasting \cite{pathak2022fourcastnet}, large-scale $\mathrm{CO_{2}}$ injection simulations \cite{Grady2022ModelParallelFN}, subsurface two-phase oil/water
flow simulations \cite{zhang2022fourier}, and so on.
\begin{figure}[tbhp]
\begin{minipage}[b]{1\linewidth}
    \centering
    \subfloat[][Example \uppercase\expandafter{\romannumeral1} (OPNO, ours)]{\includegraphics[width=.33\linewidth]{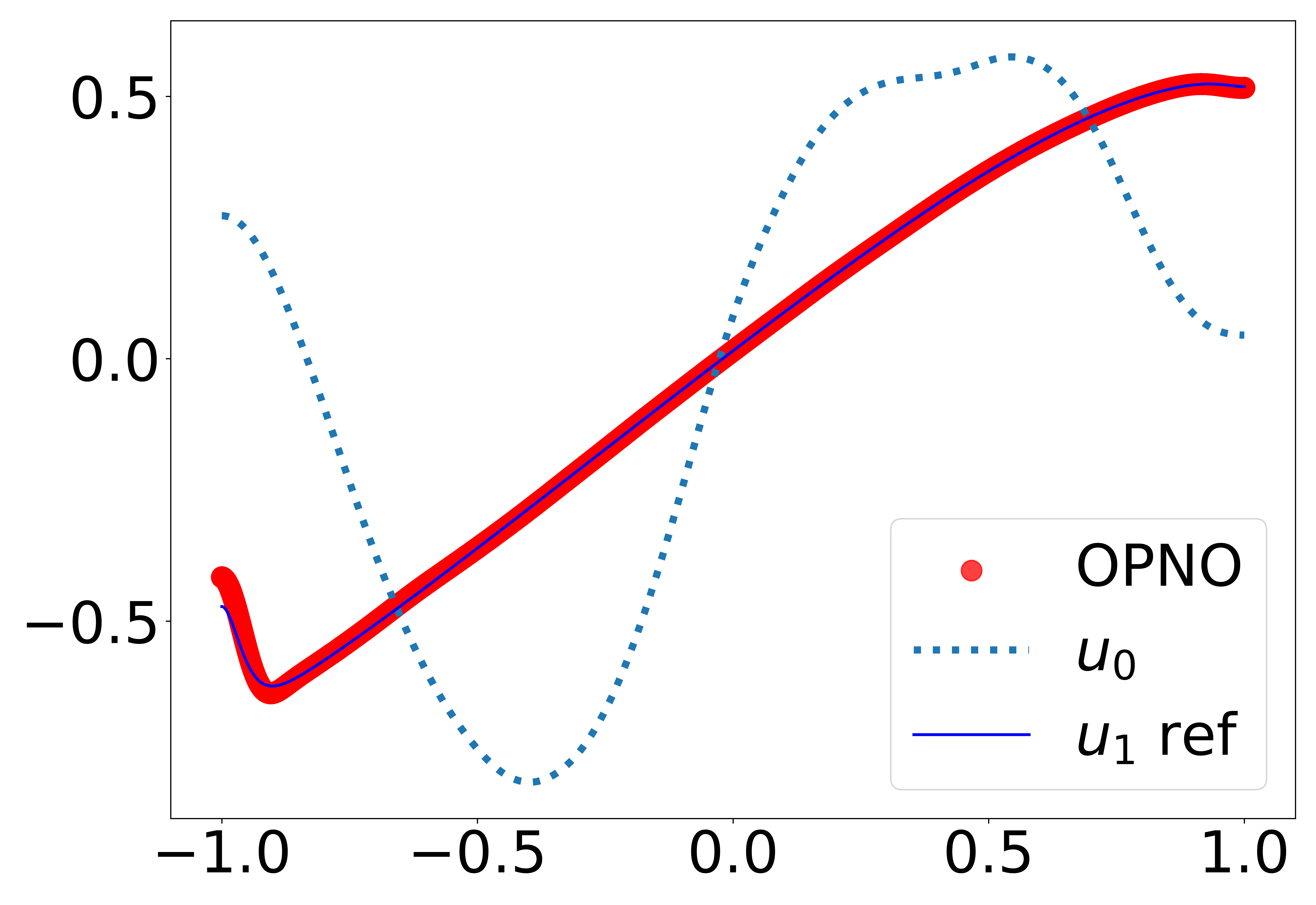}}
    \subfloat[][Example \uppercase\expandafter{\romannumeral2} (OPNO, ours)]{\includegraphics[width=.33\linewidth]{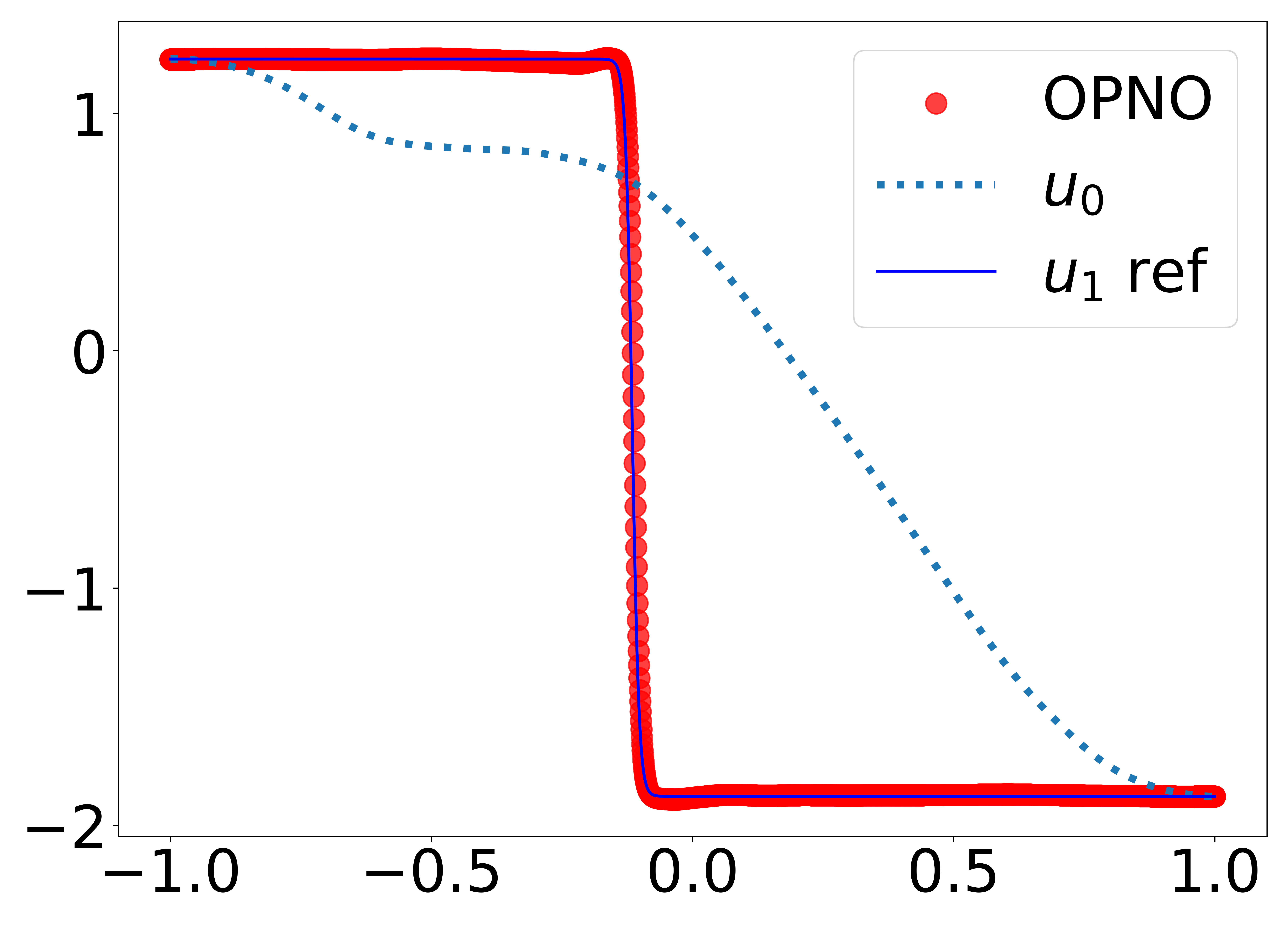}}
    \subfloat[][Example \uppercase\expandafter{\romannumeral3} (OPNO, ours)]{\includegraphics[width=.33\linewidth]{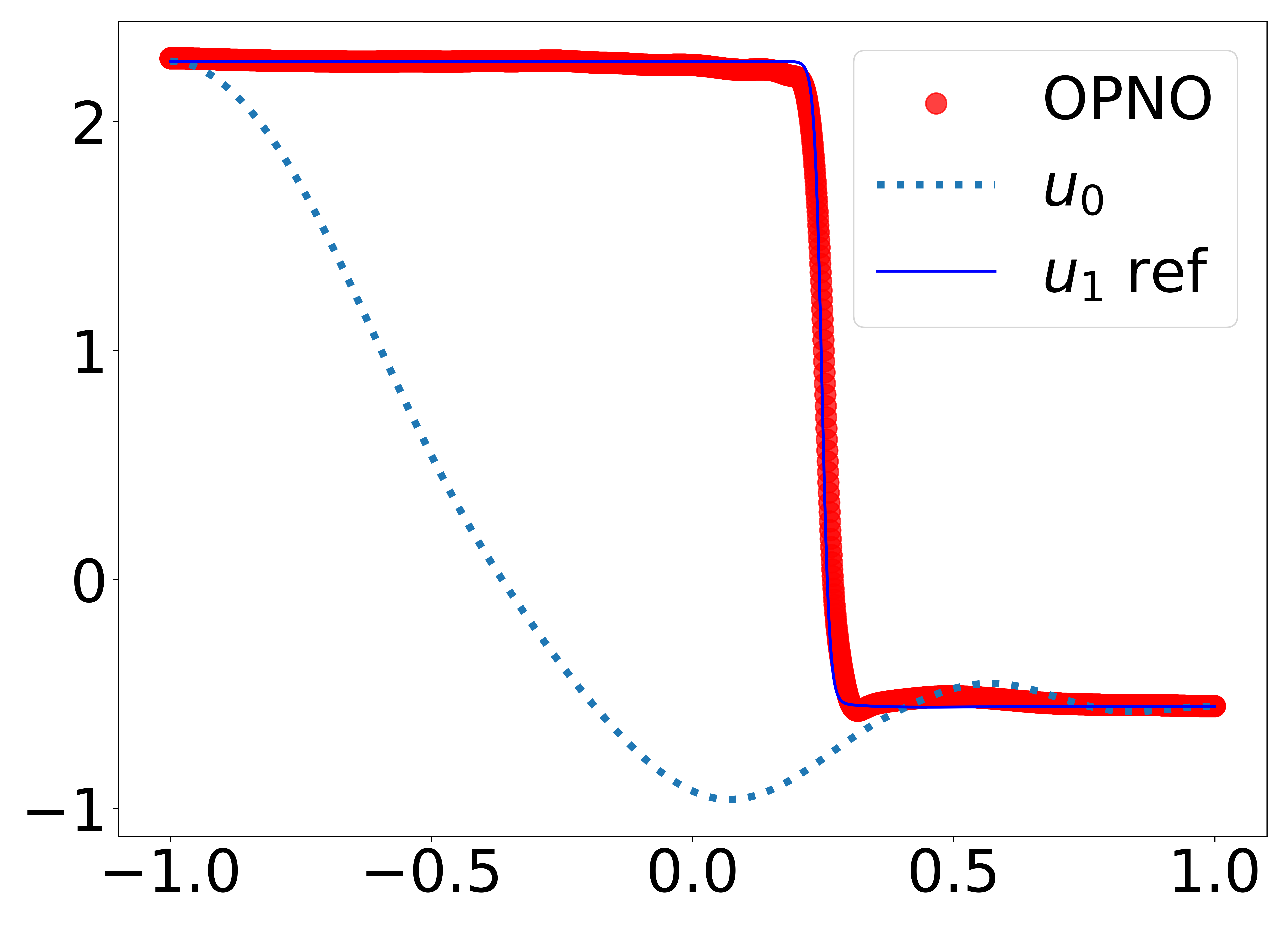}}
  \end{minipage}
\begin{minipage}[b]{1\linewidth}
    \centering
    \subfloat[][Example \uppercase\expandafter{\romannumeral1} (FNO)]{\includegraphics[width=.33\linewidth]{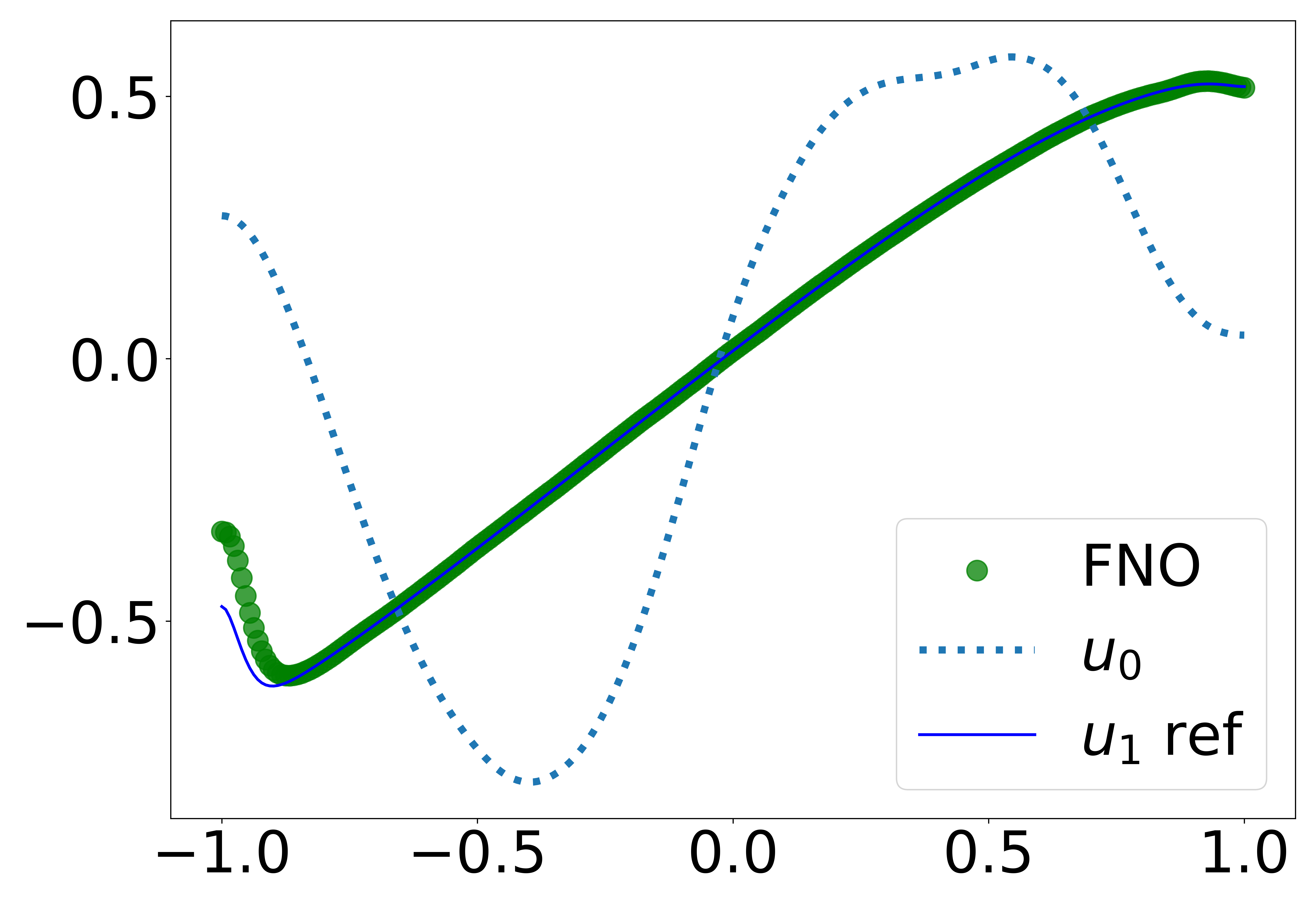}}
    \subfloat[][Example \uppercase\expandafter{\romannumeral2} (FNO)]{\includegraphics[width=.33\linewidth]{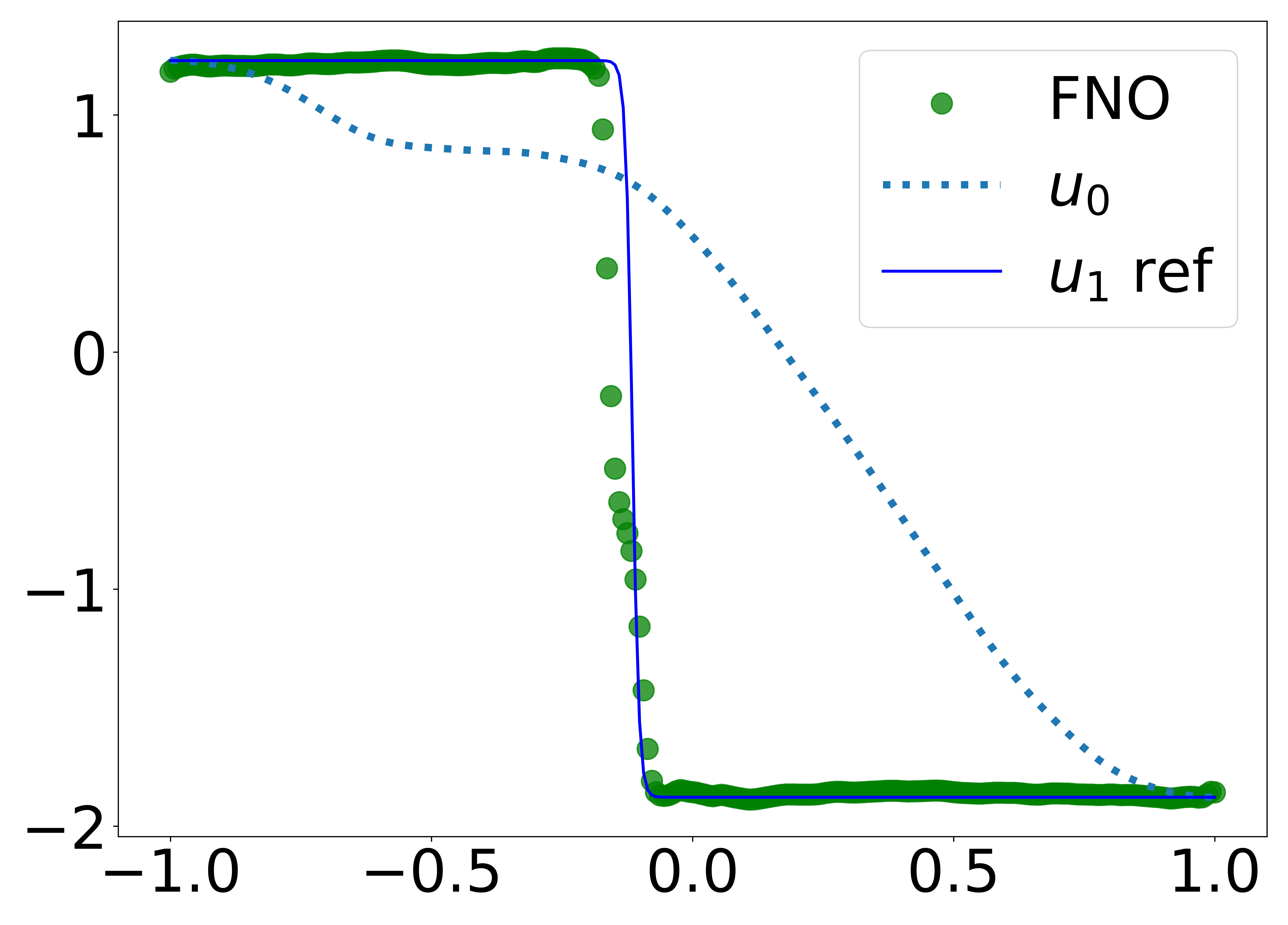}}
    \subfloat[][Example \uppercase\expandafter{\romannumeral3} (FNO)]{\includegraphics[width=.33\linewidth]{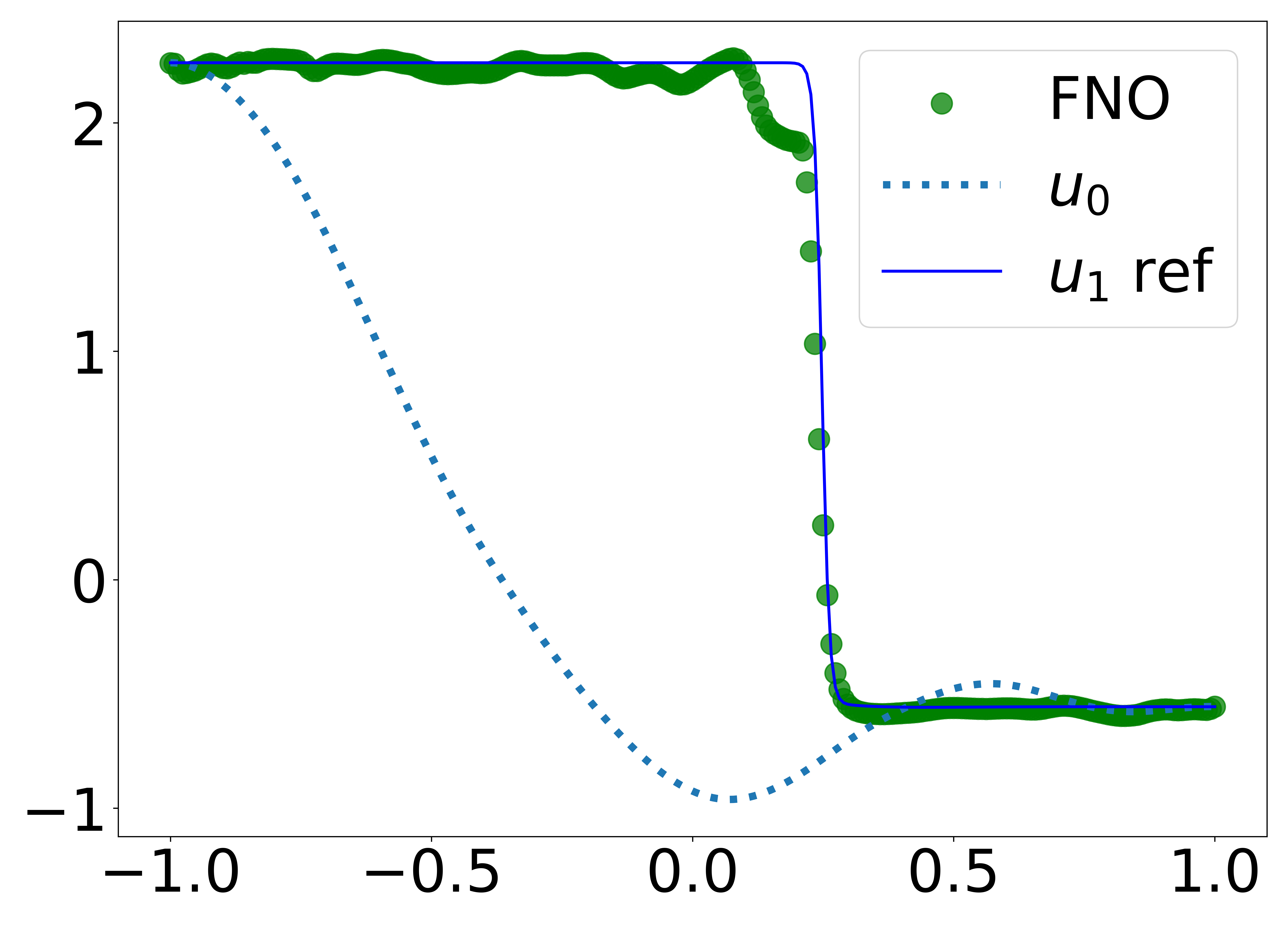}}
  \end{minipage}
  \caption{\textbf{Predictions of FNO model versus our OPNO for Burgers' equation with Neumann BCs.} The dotted and solid blue lines represent the initial condition $u_0(x)$ in the test dataset and the corresponding reference solution $u_1(x)$ when
$t=1$, respectively. The red dots in
the upper figures are the predictions of our proposed OPNO method, while the green dots in the lower figures are their counterparts produced by the FNOs. Detailed discussions are presented in \cref{exp:burgers}.}
\label{fig:neumann}
\end{figure}

While most neural operators are devoted to solving PDEs with periodic boundary conditions (BCs), limited work has been done in cases with non-periodic BCs. More precisely,
let $D \subset \mathbb R^d$ be a domain. We are interested in finding a $\theta$-learnable-parameterized neural approximation $\mathcal G_{\theta}$ for the continuous operator
$$\mathcal G:\mathcal A(D;\mathbb R^{d_a}) \to
\mathcal U(D;\mathbb R^{d_u}), \ a \mapsto u \triangleq \mathcal G(a),$$
where $a \in \mathcal A(D; \mathbb R^{d_a})$ is an input function. When $\mathcal G$ is the solution operators of PDEs, it is essential for $\mathcal G_{\theta}$ to
precisely satisfy the imposed BCs. For
example, consider the following time-dependent PDE with Neumann BCs:
\begin{equation}
\label{eq:4}
  \partial_t u(x, t) + \mathcal N(u) = 0, \ x \in D, \nonumber
\end{equation}
subject to
\begin{eqnarray}
  &u(x, 0) = u_0(x), \ x \in D, \nonumber \\ 
  &\frac{\partial u}{\partial x} (x, \cdot) = 0, \ x \in \partial D \label{eq:example_bc},
\end{eqnarray}
and denote by $S(t)$ the solution operator that evolves the initial condition $u_0$ to the solution at time $t$, i.e.,
\begin{equation}
u(x, t) = S(t) u_0, \ S(t+\delta t) = S(t)S(\delta t), \ t, \delta t \geq 0. \nonumber
\end{equation}
If $\mathcal G_{\theta}$ is a numerical approximation of $S(1)$, the output $\mathcal G_{\theta}(u_0)$ should accurately satisfy the BCs in Eq. \eqref{eq:example_bc}.

\begin{figure}[htbp]
\centering
  \centerline{\includegraphics[width=1\textwidth]{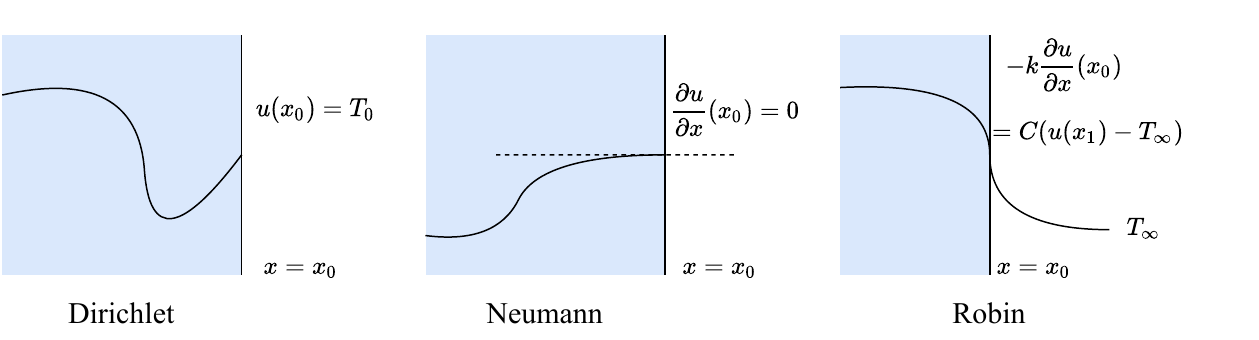}}
  \caption{\textbf{An illustration of the three types of non-periodic BCs.} The heat diffusion equation is taken as an example, for which the solution $u(x)$ represents the temperature at point $x$. While the Dirichlet BC enforces a fixed temperature on the boundary, the
Neumann BC describes an adiabatic or symmetric boundary effect by forcing the derivative to be zero, and the Robin BC, a linear combination of the two previous BCs, represents the heat
convection between the heat source and the environment on the surface.}
\label{fig:three-BCs}
\end{figure}

%Denote the unit normal on the surface by $\mathbf n$. 
The three most common types of BCs in scientific and engineering computation are the Dirichlet BC
\begin{equation}
\label{eq:5}
  u(x, \cdot) = f(x), \ x \in \partial D, \nonumber
\end{equation}
the Neumann BC
\begin{equation}
\label{eq:6}\nonumber
 \frac{\partial u}{\partial \mathbf{n}} = g(x), \ x \in \partial D,
\end{equation}
and the Robin BC
\begin{equation}
\label{eq:7}\nonumber
a u(x) + b \frac{\partial u}{\partial \mathbf{n}}(x) = h(x), \ x \in \partial D,
\end{equation}
which are known as the first-, second- and third-type boundary conditions, respectively.
\cref{fig:three-BCs} displays an example for these BCs.

Approaches for developing PINNs with hard constraints for non-periodic BCs or asymptotic limits
were introduced in \cite{lyu2022mim,zhoutao2022,liuunified,lu2021physics,jiangsong2022model}, where the neural networks needed to be retrained for different PDE parameters or initial conditions; thus, the direct generalization of these methods to neural operators seems ineffective. Besides, similar concept of preserving the properties of original system for deep-learning methods has been introduced in multiple papers by Pengzhan Jin, Yifa Tang, et al. \cite{jin2020sympnets,zhu2022approximation}. In Section \ref{sec:sol}, based on the neural network and spectral Galerkin methods \cite{shen1995efficient,shen2011spectral}, we construct a general neural framework called spectral operator learning (SOL), and
then in Section \ref{sec:opno}, we introduce an SOL method
named the orthogonal polynomial neural operator (OPNO), which generates solutions satisfying the general BCs mentioned above \textbf{up to a machine precision limit}, as well as its universal approximation
theorem. Moreover, the proposed method
possesses the following appealing properties, which are also illustrated by numerical experiments in Section \ref{sec:experiment}:

\textbf{Quasi-linear computation complexity}. The OPNO requires $O(N^d \log N)$ time complexity for $d$-dimensional problems due to the applications of the fast Chebyshev transform and the so-called fast compacting
trasform.

%In numerical experiments, it cost approximately
%$1.5$ times as much time as the FNO when under the same settings, mainly because in the code of current version, the discrete cosine transform is implemented by employing unmodified FFT after
%extending the input vectors to $2N$-length for simplicity, as in \cite{trefethen2000spectral}.

\textbf{Efficient spatial differentiation computation}. Since the output of the OPNO can be viewed as a polynomial function, its derivatives of arbitrary orders can be accurately computed within
$O(N\log N)$ operations by the numerical method named
differentiation in the frequency
space (see \cref{app:chebyshev}). This approach avoids the vanishing gradient problem in the automatic differentiation of neural networks. In addition, readers may find that this property leads to a zero-shot
learning methodology for solving
PDEs, and this will be further investigated in our future work.

\textbf{No overfitting during training}. We have observed that as two typical examples of SOL, both the FNO and OPNO are well self-regularized when solving all the PDEs this paper involves, which means that their test errors
do not obviously increase after extremely long-term training, even on a tiny training dataset. It is an extremely favorable feature for a deep learning method and suggests that an ultra-precise SOL
method for specific PDEs is practically feasible.

\textbf{Quasi-spectral accuracy}.
The spectral structure of the SOL approach leads to a behavior that is analogous to the spectral accuracy of spectral methods, which means that under the assumption of smoothness, a model trained on a coarse mesh can be directly
applied to generate solutions on fine meshes without loss of numerical accuracy, and we presume that this is the reason for the observation of the ``\textbf{resolution-invariant}''\cite{li2021fourier,gupta2021multiwavelet,cao2021choose,tripura2022wavelet} properties of some FNO-like methods.

Both of the code and dataset for the paper are available at \\ \url{https://github.com/liu-ziyuan-math/spectral\_operator\_learning}.
\section{Spectral Operator Learning}\label{sec:sol}

A general spectral neural operator $\mathcal G_{\theta} : \\ \mathcal A(\Omega;\mathbb R^{d_a}) \rightarrow \mathcal U(\Omega;\mathbb R^{d_u}), a \mapsto G_{\theta}(a)$ is a mapping of the form
\begin{equation}
  \mathcal G_{\theta}(a) =\mathcal Q \circ \mathcal L^{(L)}  \circ \sigma \circ \mathcal L^{(L-1)} \circ ... \circ \sigma \circ \mathcal L^{(1)} \circ \sigma \circ \mathcal P(a), \nonumber
\end{equation}
where $\sigma$ is a non-polynomial activation function; $\mathcal Q:\mathcal A(D;\mathbb R^{d_a}) \rightarrow \mathcal U(D; \mathbb R^{d_v})$ and $\mathcal P:\mathcal A(D;\mathbb R^{d_v}) \rightarrow \mathcal U(D;
\mathbb R^{d_u})$ are some simple structured neural networks; and $\mathcal L^{(l)}$ is a linear spectral operator layer of the form
\begin{equation}
\label{eq:3}
 \mathcal L^{(l)}(v) = \mathcal T^{-1}( A_l \cdot \mathcal T (v) ).
\end{equation}
\begin{figure}[tbhp]
\centering
\centerline{\includegraphics[width=9cm]{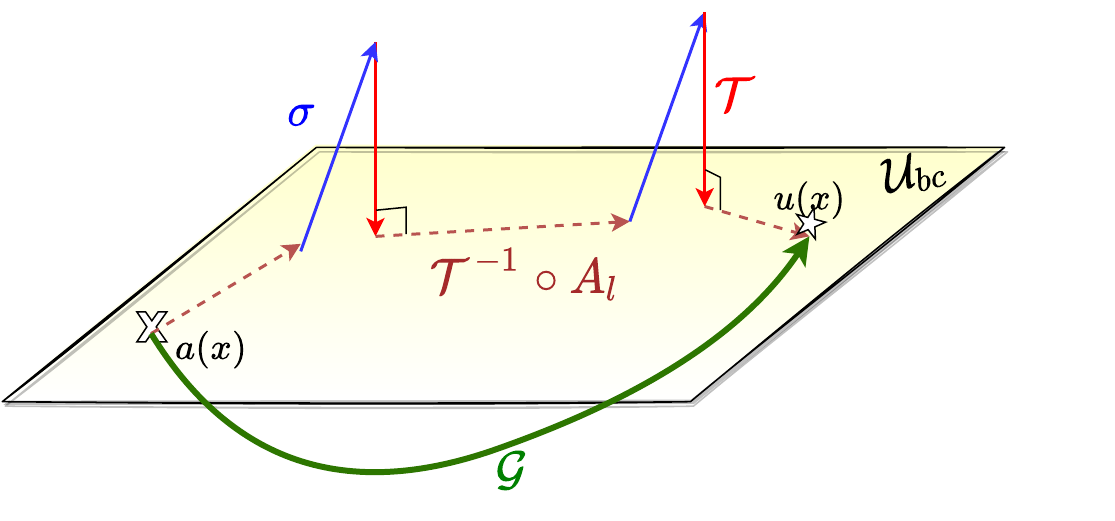}}
\caption{\textbf{Illustration of the architecture of SOL from the function space perspective}. Noting that $\mathcal T$ is the projection to the function space of underlying basis that
  satisfy specific BCs, while $\mathcal T^{-1} \circ A_l$ is a linear transform on the space.}
\label{fig:sol}
\end{figure}
In Eq. \cref{eq:3}, $\mathcal T$ is a specific transformation operator that decomposes functions into
the frequency domain of the corresponding basis, for which the inverse operator is denoted by
$\mathcal T^{-1}$, and $A_l \in \mathbb R^{d_v \times dv}$ is a (learnable) parameterized linear transformation of the frequency domain. The idea of such an architecture was first introduced by Zongyi Li, Nikola Kovachki, et al. in the FNO \cite{li2021fourier}, where $\mathcal T = \mathcal F$ (the Fourier
transform) was considered, and $A_l$ was fixed as a diagonal matrix:
\begin{equation}
\label{eq:8}
(\mathcal K(\phi) v) (x) \triangleq \mathcal F^{-1} \big( \mathrm{diag}(\Lambda) \cdot (\mathcal F v) \big) (x).
\end{equation}
Its significant associations with the numerical methods and the importance of choosing a correct basis, however, have not been fully discussed, which we find lead to a general framework of learning-based methods.

\begin{figure}[htbp]
\centering
    \subfloat[][Fourier method]{\includegraphics[width=.25\linewidth]{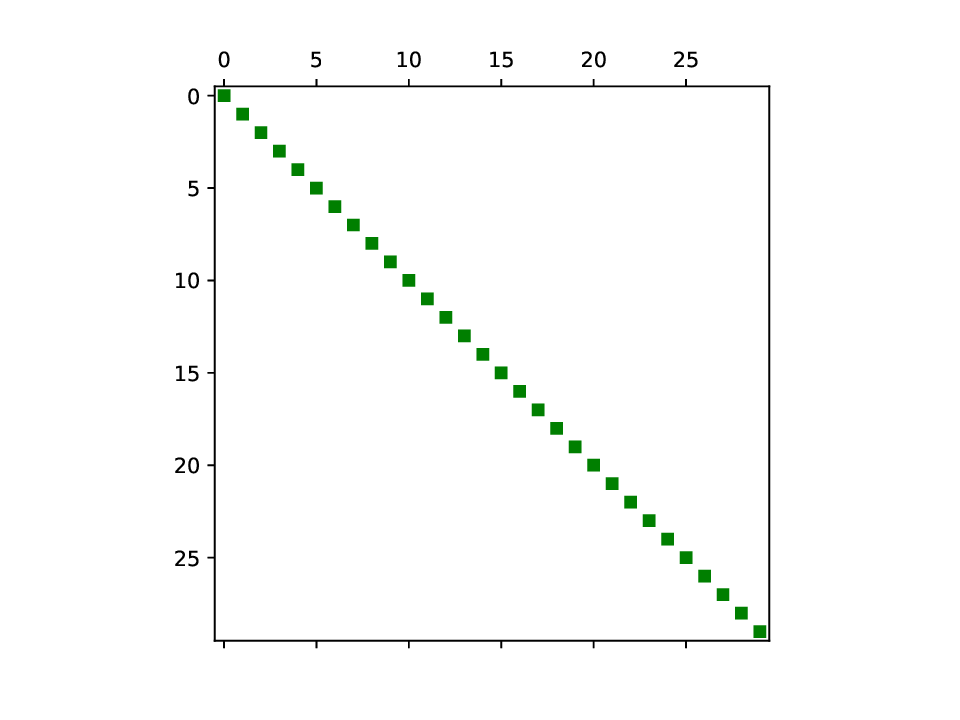}}
    \subfloat[][FNO]{\includegraphics[width=.25\linewidth]{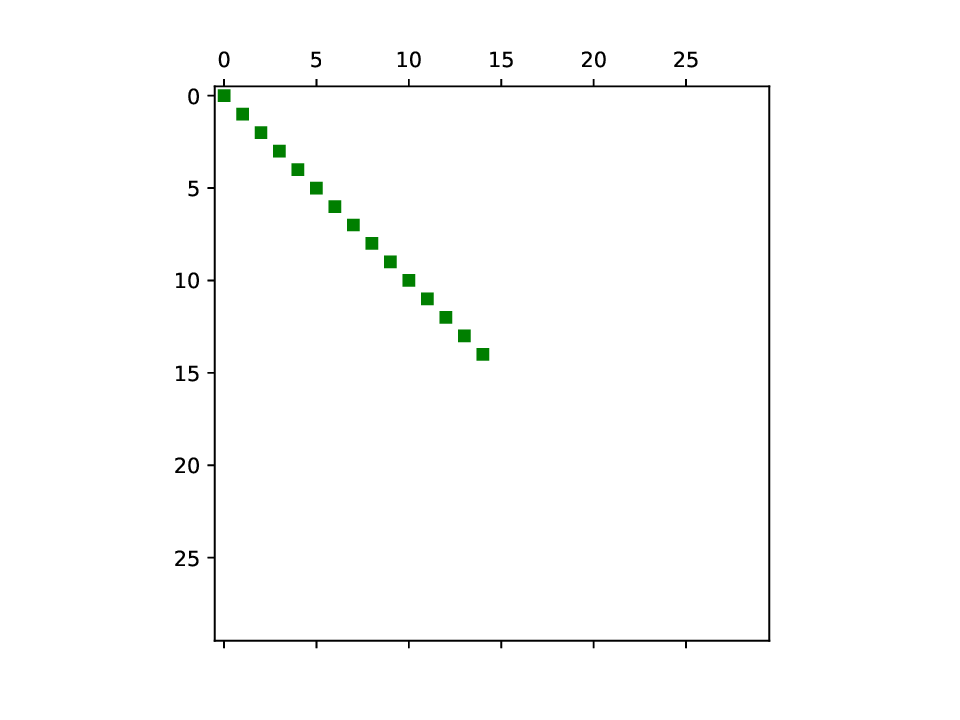}}
    \subfloat[][(Legendre) Garlerkin method]{\includegraphics[width=.25\linewidth]{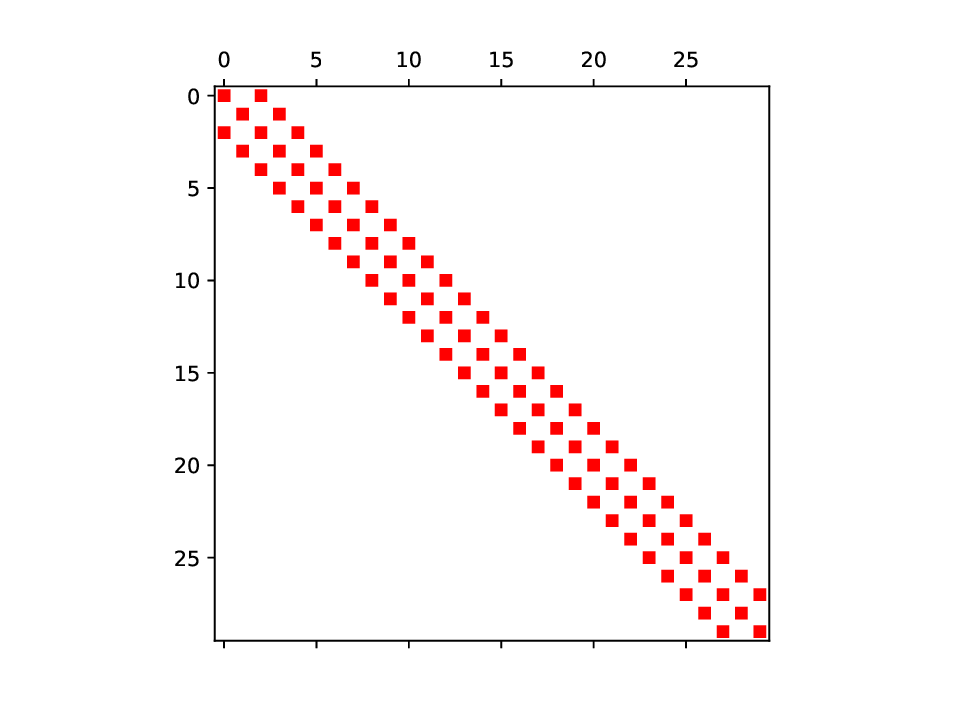}}
    \subfloat[][OPNO (ours)]{\includegraphics[width=.25\linewidth]{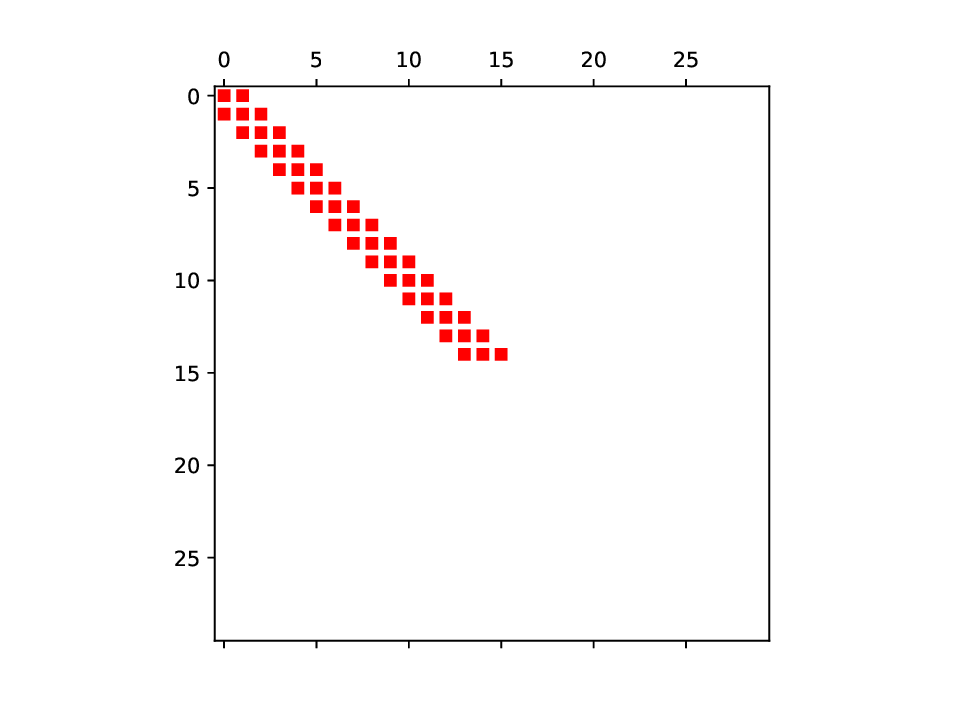}}
    \caption{\textbf{Illustration for the sparse (spectral) linear systems} of the pseudospectral Fourier method (a), single-layer FNO (b), spectral Galerkin method (c) and single-layer OPNO (d), when solving the example boundary value problems
      $\Delta u + \alpha u = 0, \ \alpha \in \mathbb R$.}
   \label{fig:linear-system}
\end{figure}

 \cref{fig:sol} illustrates a sketch map for the SOL scheme. One can determine that the architecture of SOL alternately transforms the function \textbf{linearly} in the
solution (frequency) space by a spectral operator and maps the function \textbf{nonlinearly} onto
the physical space
by an activation function.
Readers may immediately discover the resemblance between SOL and the so-called pseudospectral techniques in numerical method, which numerically solves the \textbf{linear} part of a PDE with the
spectral method in the frequency domain, while the \textbf{nonlinear} terms are solved in the physical space. Moreover, the spectral linear systems of SOL share the same pattern with their corresponding
spectral methods (see  \cref{fig:linear-system}).

\begin{figure}[htbp]
\centering
\subfloat[][$y=tanh(4x)$]{\includegraphics[width=.45\linewidth]{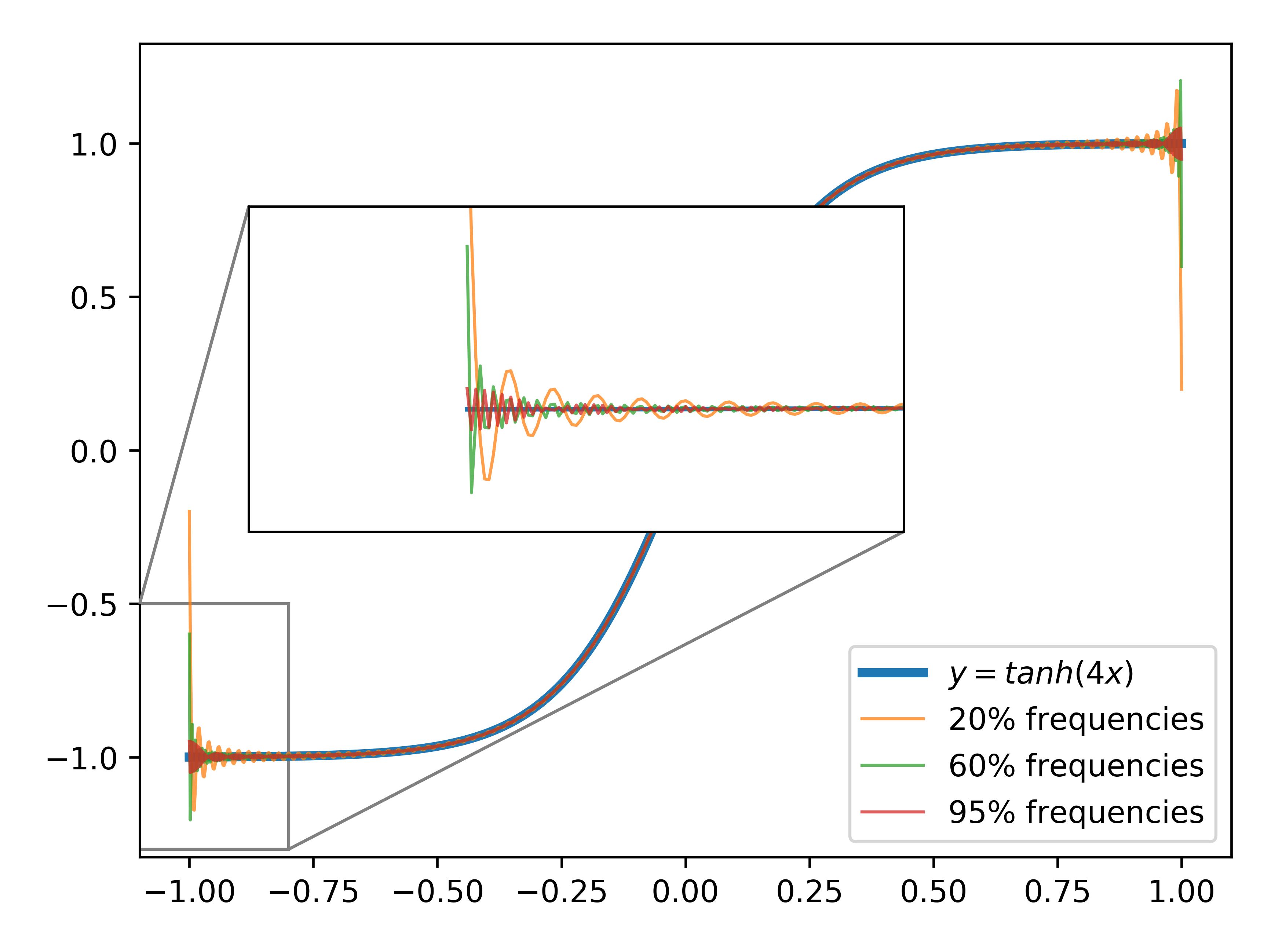}}
\subfloat[][$y=x$]{\includegraphics[width=.45\linewidth]{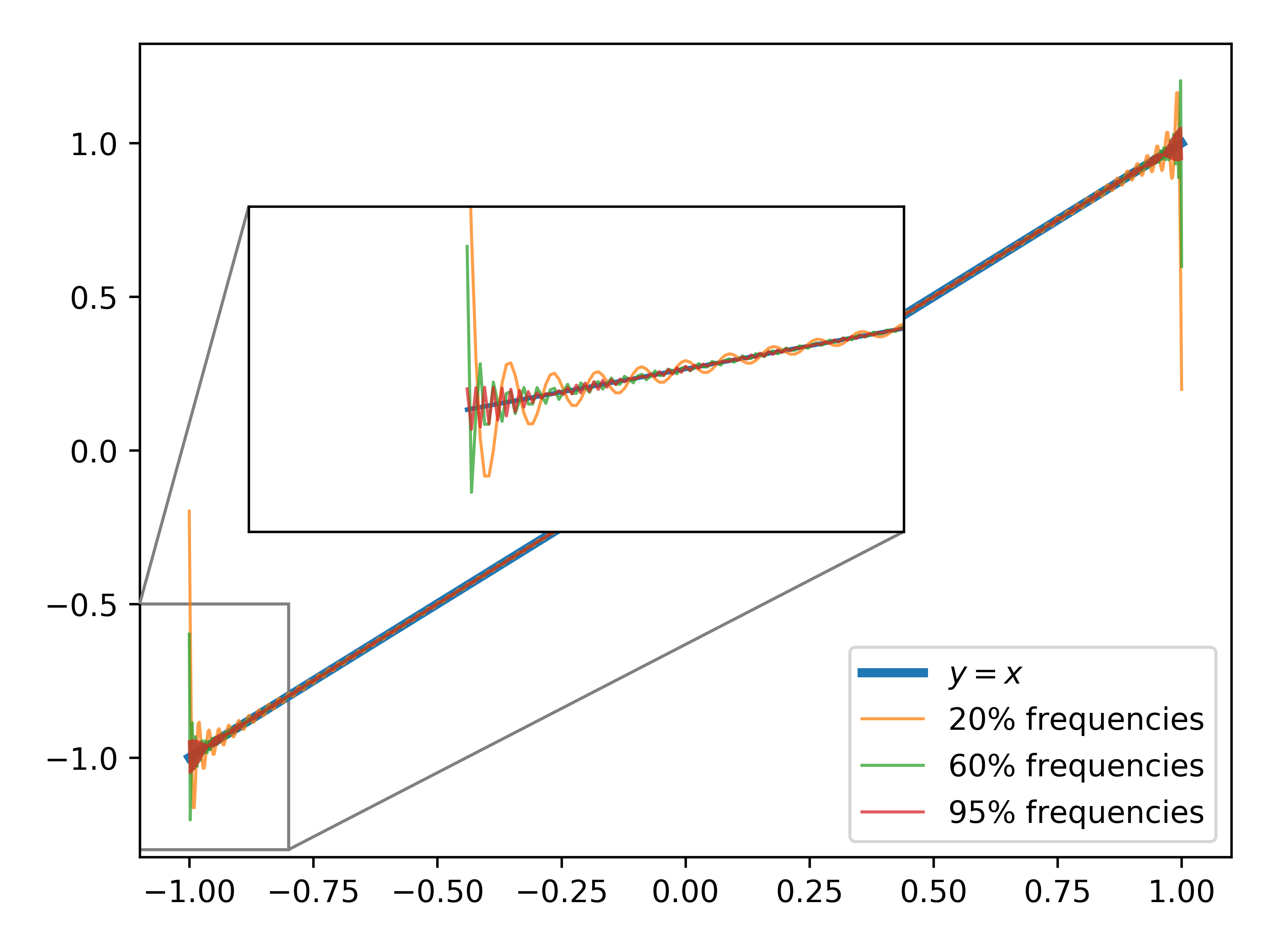}}
\caption{\textbf{The truncated Fourier expansions of non-periodic functions suffer from Gibbs phenomona}. The number of discretization points is set as $N = 1000$.}
\label{fig:gibbs}
\end{figure}

While the FNO is now widely applied in many problems with large scales and high complexity, most of the problems it has effectively solved are periodic,
and we find that the performance of the FNO drops once the assumption of periodicity is not satisfied. For instance, \cite{ong2022iae} indicated that the FNO is not suitable for solving scattering problems with Sommerfeld radiation
conditions. In addition, since the FNO is built on the truncated basis of trigonometric polynomials, the accuracy and computational stability may also
be troubled with the so-called Gibbs phenomenon, which means that when applying Fourier methods to a non-periodic problem, spurious high-frequency oscillations will be generated near the boundaries, and the global convergence rate will
be severely reduced (see  \cref{fig:gibbs}).
It is common sense that for these kinds of PDEs, spectral methods with other bases, such as Chebyshev, Legendre and Hermite polynomials, should be applied rather than the Fourier spectral
method \cite{shen2011spectral}.

Therefore, based on the spectral methods of orthogonal polynomials, especially of their compact combinations that will be introduced in \cref{sec:c_p}, we develop the OPNO to
efficiently and reliably solve non-periodic PDEs that are subject to the Dirichlet, Neumann, and Robin BCs. We focus on Chebyshev polynomials throughout this paper. The methodology,
however, is also feasible for neural operators with other bases such as Legendre, Jacobi, Laguerre, Hermite polynomials, or spherical harmonic functions for specific problems, e.g., those on the cylindrical, spherical, or unbounded domain. However, such issues will not be addressed here. Furthermore, hybrid bases are useful for solving separable multi-dimensional problems or coupled equations with different BCs under the SOL framework.

\section{The OPNO and its fast algorithm}\label{sec:opno}
\begin{figure}[tbhp]
\centering
  \centerline{\includegraphics[width=1\textwidth]{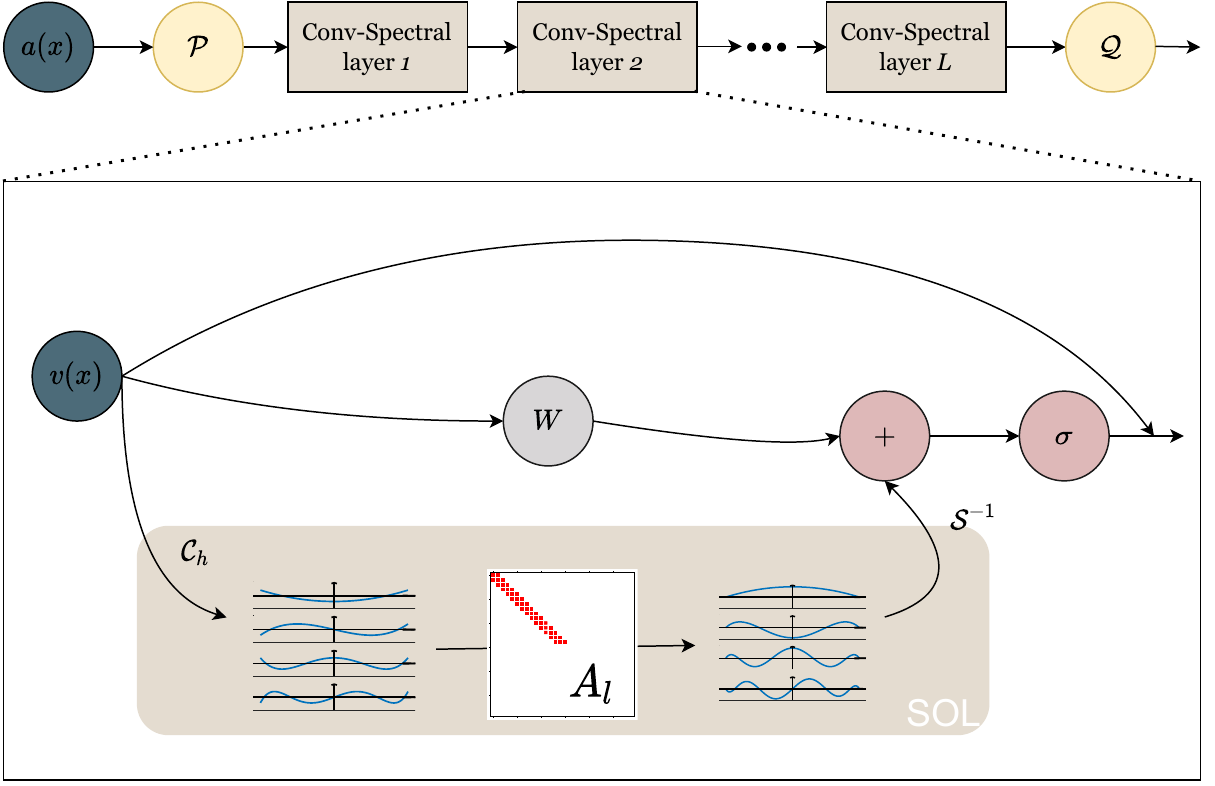}}
  \caption{The architecture of the OPNO layers.}
  \label{fig:OPNO}
\end{figure}
Without any loss of generalization, the interval $I = [-1, 1]$ is considered since the proposed method can be easily generalized to any interval $[a, b] \subset \mathbb R$ by a linear map such that
\begin{equation}
t = x(b-a)/2 + (a+b)/2, \ t \in [a, b] \nonumber
\end{equation}
or to $d$-dimensional cases.
Based on the compact combination basis of Chebyshev polynomials $\left\{ \phi_k(x) \right\}$ that Jie Shen, Tao Tang and Li-lian Wang introduced in \cite{shen2011spectral}, we construct the fast Shen transform
$\mathcal S$ that maps an arbitrary function $u \in C(I)$ to its expansion coefficients on the compact combination basis.
To fulfill this goal, the Shen transform is split into the composition $\mathcal S = \mathcal C_p \circ \mathcal
C_h$: the fast Chebyshev transform $\mathcal C_h$ that maps the function into its Chebyshev coefficients is taken first ($x\rightarrow T_k$), and then the fast compacting transform $\mathcal C_p$ mapping the Chebyshev
coefficients to the coefficients of compact combinations $\phi_k$ ($T_k\rightarrow\phi_k$) follows.

\subsection{Brief introduction to the fast Chebyshev transform $\mathcal C_h$}
Chebyshev polynomials (of the first kind) $T_n(x)$ are a seires of orthogonal polynomials given by the following three-term recurrence relation:
\begin{equation}
\label{eq:11}
T_{n+1}(x) = 2 x T_n(x) - T_{n-1}(x), \ n \geq 1, \ x \in I.
\end{equation}
To list a few, $T_0(x)=1, \ T_1(x) = x, \ T_2(x)=2x^2-1$, and $T_3(x) = 4x^3-3x$. Denote $x_j = - \cos \frac{\pi j }{N}, \ 0 \leq j \leq N$ as the Chebyshev-Gauss-Lobatto (CGL) points, and
let $I_N u$ be the
Lagrange interpolation polynomial relative to these CGL points; then, any function $u \in C(I)$ has a unique decomposition on the basis of $T_k$, namely,
\begin{equation}
(\Pi_N u)(x) = (I_Nu)(x) \triangleq \sum_{n=0}^N \hat u_n T_n(x), \nonumber
\end{equation}
where $\Pi_N u$ denotes the Chebyshev decomposition transform and the Chebyshev coefficients $\left\{ \hat u_n \right\}$ are determined by the (forward) discrete Chebyshev transform $\mathcal{C}_{h}$. The mathematical forms of $\mathcal C_h$ and its inverse are
listed below.
\begin{equation}
\begin{split}
  &\hat u_n = \frac{2}{\tilde c_n N} \sum_{j=0}^N \frac{1}{\tilde c_j} u(x_j) \cos \frac{nj\pi}{N}, 0 \leq n \leq N \ \text{(forward transform)}, \\
  &u(x_j) = \sum_{n=0}^N \hat u_n \cos \frac{n j \pi}{N}, 0 \leq j \leq N \ \text{(backward transform)}, \nonumber
\end{split}
\end{equation}
where $\tilde c_0 = \tilde c_N = 2$ and $\tilde c_j = 1$ for $j=1, 2, ..., N-1$.

Moreover, the discrete Chebyshev transform and its inverse can be efficiently computed in $O(N log N)$ operations via the fast cosine transform \cite{trefethen2000spectral,wen1977fast}. A more detailed introduction is given in \cref{app:chebyshev}. 

\subsection{Compacting transform $\mathcal C_p$}\label{sec:c_p}
As \cite{shen2011spectral} introduced, the compact combination of Chebyshev polynomials $\left\{ \phi_k(x) \right\}$ is a basis of polynomails that automatically satisfies specific BCs; this combination has the following form
\begin{equation}
\phi_i(x) = T_k(x) + a_k T_{k+1}(x) + b_kT_{k+2}(x). \nonumber
\end{equation}
More precisely, consider the general (Robin) BCs:
  \begin{equation}
a_-u(-1) + b_- u'(-1) = 0, \ a_+ u(1) + b_+ u'(1) = 0.  \nonumber
\end{equation}
Then, there exists a unique set of $\left\{ a_k, b_k \right\}$ such that $\phi_k(x)$ satisfies such conditions for all $k > 0$.
It is remarkable that:
\begin{itemize}
\item For $a_{\pm} = 1$, $b_{\pm} = 0$ (Dirichlet BCs), we have $\phi_k(x) = T_k(x) - T_{k+2}(x)$.
\item For $a_{\pm} = 0$, $b_{\pm} = 1$ (Neumann BCs), we have  $\phi_k(x) = T_k(x) - \frac{k^2}{(k+2)^2}T_{k+2}(x)$.
\end{itemize}

The compacting transform $\mathcal{C}_{p}$ maps the Chebyshev coefficients to the expansion coefficients on the $\phi_k$ of the same function. Notably, the compacting transform and its inverse are linear and can be carried out by
recursion with $O(N)$ operations. For example, let $f = \sum_{j=0}^{N} \alpha_j T_j = \sum_{j=0}^{N-2} \beta_j \phi_j $, and consider the Dirichlet BCs such that $\phi_k = T_k - T_{k+2}$;
then, we have the following:

\textbf{forward compacting transform}
\begin{equation}\label{eq:forward_cp}
\beta_j =\left\{\begin{aligned}
  &\alpha_j, \ j = 0, 1, \\
  &\beta_{j-2} + \alpha_j, \ 2 \leq j \leq N-4, \\
  &-\alpha_{j+2}, \ j = N-3, N-2, \\
\end{aligned}\right.
\end{equation}

and \textbf{backward compacting transform}
\begin{equation}
\alpha_j =\left\{\begin{aligned}
    &\beta_j, \ j = 0, 1, \\
    &\beta_j - \beta_{j-2}, \ 2 \leq j \leq N-2, \\
    &- \beta_{j-2}, \ j = N-1, N.
\end{aligned}\right. \nonumber
\end{equation}
The forward and backward compacting transforms for cases with Neumann BCs are provided in \cref{app:cp-neumann}.

\subsubsection{Fast compacting transform}
\begin{figure}[tbhp]
  \centerline{\includegraphics[width=0.7\textwidth]{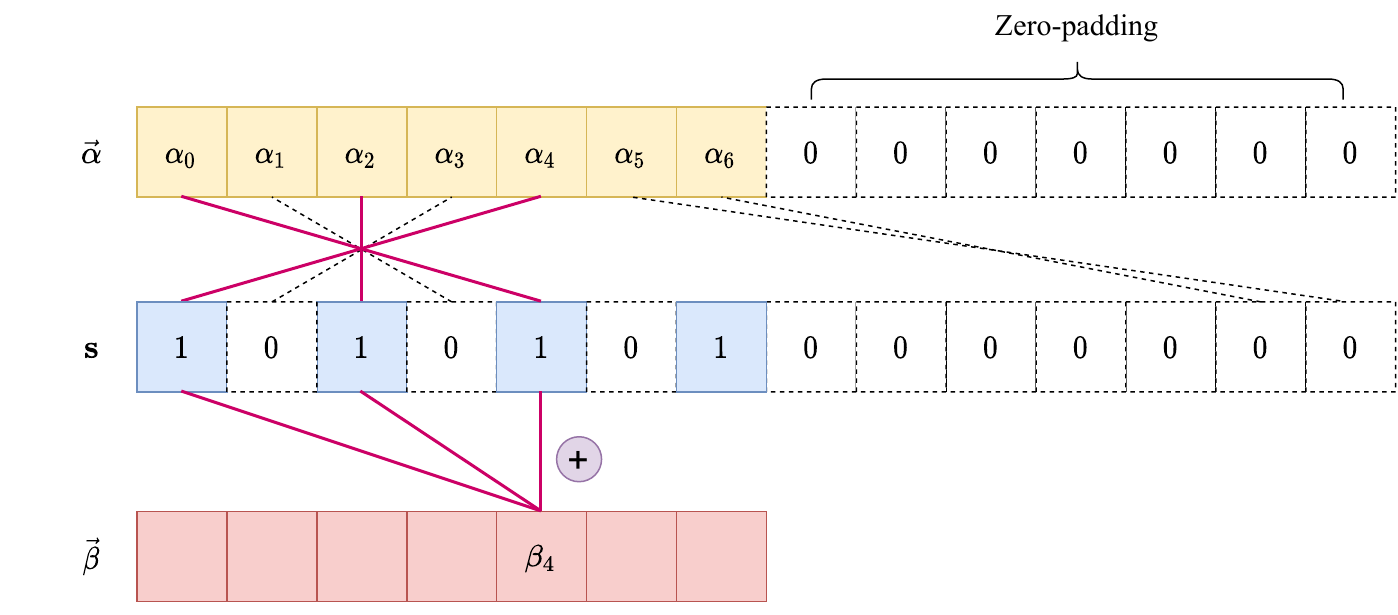}}
  \caption{\textbf{Schematic diagram of the linear convolutions of fast compacting transform} $\mathcal C_p$ when $N=6$, so that $\mathbf{s}=(1,0,1,0,1,0,1)^T$ and $\beta_4=\alpha_4+\alpha_2+\alpha_0$.}
  \label{fig:linear-conv}
\end{figure}

Although the native forward compacting transform theoretically costs $O(N)$ floating point operations (FLOPs), the coefficients need to be computed by recursion; thus, it is not suitable for computation on GPUs. Consequently, we
propose a parallelable fast compacting transform. Suppose $N \geq 3$, and let $\mathbf{s}$ be an $(N+1)$-dimensional vector with elements
\begin{equation}
s_j = \left\{\begin{aligned}
    &0, \ j = 1, 3, 5, ..., 2\lceil N/2 \rceil - 1, \\
    &1, \ j = 0, 2, 4, ..., 2\lfloor N/2 \rfloor.
\end{aligned}\right. \nonumber
\end{equation}
We usually ignore the special cases of $j \geq N-3$ in Eq. \cref{eq:forward_cp} for simplicity or computational stability when the function $f = \sum_{j=0}^{N} \alpha_j T_j$ satisfies the
corresponding BCs. Then one can see that Eq. \cref{eq:forward_cp} is equivalent to the \textbf{linear convolution} of $\boldsymbol{\alpha}$ and $\mathbf{s}$, i.e.,
$$
\beta_j = \sum\limits_{k=0}^{\lfloor j/2 \rfloor} \alpha_{j-2k} = \sum\limits_{k=0}^{j} \alpha_{j-k} s_{k}.
$$
The fast algorithm for discrete linear convolution is alright highly developed and can
be quickly given by the FFT algorithm, namely
\begin{equation}
 {\boldsymbol{\beta}} = \mathcal F^{-1} (\mathcal F(\bar{\mathbf s}) \cdot \mathcal F(\bar{\boldsymbol{\alpha}}) ) ,\nonumber
\end{equation}
where the $\bar \bullet$ symbol represents zero-padding the vector to a length of $2N+2$ at its end (see  \cref{fig:linear-conv}).
The fast compacting transform exhibits approximately a threefold increase in speed in the numerical experiments.

\subsection{Structure of the OPNO}
We finally come to the structure of the OPNO. The vanilla spectral kernel layer is first given by substituting $\mathcal T$ in Eq. \cref{eq:3} with the Shen transform $\mathcal S \triangleq \mathcal{C}_p \circ \mathcal{C}_{h}$, i.e.,
\begin{equation}
  \tilde{\mathcal L}_{ker}^{(l)}(v) = \mathcal S^{-1}( A_l \cdot \mathcal S (v) ). \nonumber
\end{equation}
In addition, as the linear system of the spectral Galerkin method implies (see  \cref{fig:linear-system}), $A_l$ is a quasi-diagonal matrix with a bandwidth of $w$, and truncated to $k$ modes to reduce
the number of model parameters and focus on the low-frequency behavior of
the target operator. We find that $w=3$ suffices in all the numerical examples. In addition, out of consideration for efficiency, since $\mathcal
C_{p}$ is such a simple linear transform that can be easily learned, it is omitted in all forward Shen transform steps except in the final projection layer $\mathcal Q$. In total, a further 2X speed
optimization is witnessed without noticeable accuracy loss. Therefore, the spectral kernal layer we implement is given by
\begin{equation}
\mathcal L_{ker}^{(l)}(v) = \mathcal S^{-1}( A_l \cdot \mathcal C_h (v) ). \nonumber
\end{equation}

On top of that, there are some details worth mentioning. As we have emphasized before, it is significant for the outputs to accurately satisfy the BC, so the output of the OPNO is again projected onto the space of
$\left\{ \phi_k \right\}$ via the projection layer $\mathcal Q$; namely,
\begin{equation}
\label{eq:24}
\mathcal Q(v) = \mathcal S^{-1} [\mathbf{1}_{k \leq N} \mathcal S (Q (v))],
\end{equation}
where $Q$ is usually a shallow neural network. In fact, to prove the
convergence properties of the FNO for periodic problems, \cite{kovachki2021universal} has introduced the same projection to FNOs on the space of trigonometric functions rather than the native $\sigma$-activated FNOs. We assume
that the native FNO sacrifices the inherent periodic property for generalizations to non-periodic problems. As the output of the OPNO is a finite summation of underlying basis functions, the following BC
satisfaction theorem holds.
\begin{thm}[BC satisfaction]
For any OPNO $\mathcal G_{\theta}$ with the projection operator $\mathcal Q$ of the form Eq. \cref{eq:24}, and any $a \in H^s(I)$, the output $u=\mathcal G_{\theta}(a)$
satisfies the BCs of the underlying basis of $\mathcal G_{\theta}$.
\end{thm}

Moreover, to capture the high-frequency part of the target operator, just as the FNO deployed, we also assemble a joint ``Conv-Spectral'' layer by adding an auxiliary linear neural layer $W_l$ parallel to each $\mathcal L_{ker}^{(l)}$
layer. From the perspective of spaital dimension, such linear neural layers are Convolutional Neural Networks (CNN) with a $1\times 1$ kernel; remark that they do not break the SOL-structure in  \cref{fig:sol}. The
joint Conv-Spectral layers are the backbone of our model, and readers may also see a sketch map shown in  \cref{fig:OPNO}.

Finally, we find that the skipping-connection skill of ResNet may slightly improve the resulting performance. To prevent the lifting operator $\mathcal P$ from interrupting the skipping structure, $\mathcal P$ is
implemeted as a dense block of the input function and a kernel layer.

In summary, the architecture of the OPNO layer is demonstrated as follows,
\begin{equation}
\label{eq:25}
 v^{(l+1)} = v^{(l)} + \sigma(W_l v^{(l)} + b_l + \mathcal S^{-1}( A_l \cdot \mathcal C_h (v^{(l)}) ) ).
\end{equation}

The OPNO also satisfies the following universal approximation theorem, and the associated proof is provided in \cref{sec:app-uat}.
\begin{thm}[Universal Approximation Theorem]\label{thm:uat}
Let $s\geq 0$ and let $\mathcal G: H^s(I^d; \mathbb R^{d_a}) \to L^2(I^d; \mathbb R^{d_u})$ be a continuous operator. Consider $K \subset H^s(I^d; \mathbb R^{d_a})$ a compact subset; then, $\forall
\epsilon > 0$, there exists an OPNO $\mathcal L: H^s(I^d; \mathbb R^{d_a}) \to L^2(I^d; \mathbb R^{d_u})$ such that
\begin{equation}
\label{eq:55}
\left\| \mathcal L (a) - \mathcal G (a)\right\|_{L^2} < \epsilon,  \ \forall a \in K. \nonumber
\end{equation}
\end{thm}

\section{Numerical Experiments}\label{sec:experiment}
To verify the accuracy and computational features (quasi-linear computational complexity, self-regularization, resolution-invariant error, etc.) of the OPNO, we compare it with four popular neural operator models and the finite difference method (FDM) in terms of solving four PDEs: (1) Burgers' equation with Neumann BCs; (2) the heat diffusion equation with Robin BCs; (3)
the heat diffusion equation with inhomogeneous Dirichlet BCs; and (4) the 2D Burgers' equation. The accuracy of the models is measured by the average relative $L^2$ norm error between the predicted solution and the reference solutions, which are
given by Chebyshev spectral methods, as well as the $L^{\infty}$ norm error on the corresponding BCs. The four baseline deep learning models are listed below.
\begin{enumerate}
  \item \textbf{FNO} is a state-of-the-art neural operator for solving parametric PDEs \cite{li2021fourier}. The original paper also presented a series of methods for generating supervised training datasets for
the learning-based PDE solvers, which are now cited as classic datasets in multiple papers.
\item \textbf{Galerkin transformer} (GT) \cite{cao2021choose}, a self-attention-based neural operator that introduces a novel linear and softmax-free attention mechanism, has wide applications in both PDE solving and pattern recognition.
\item \textbf{IAE-Net} \cite{ong2022iae}, a neural operator that consists of an encoder-decoder structure and integral transforms, has surpassed various models in solving PDEs and performing signal/image processing in experiments,
especially those with non-periodic BCs such as radiation BCs.
\item \textbf{POD-DeepONet} (POD-DO) \cite{lu2022comprehensive} is an improved DeepONet model based on the proper orthogonal decomposition (POD) basis and is used for comparison with FNO in the paper. The DeepONet model is
  acknowledged for its computational efficiency and geometric flexibility. It also has a solid theoretical foundation due to its
streamlined branch-trunk structures.
\end{enumerate}

To eliminate the impact of different
neural network training techniques, the hyperparameters of different models are selected to be the same, or consistent with those in their original papers, whereas the number of training epochs and decay step size of the learning rate simply increase to $5000$ and
$500$, respectively, to ensure that all models are well trained (see \cref{tab:models-para}). The learnable parameters of all models are in double-precision format ('float64' or 'complex128'). % We have observe that ...

%It's important for the dataset that ...[??1]

\begin{table}[htbp]
  \footnotesize
  \caption{\textbf{Default settings of the deep learning methods}: batch sizes (bs), numbers of training epochs, training optimizers, spectral modes, spectral channel widths, numbers of Conv-Spectral layers ($L$), types of interpolation grids, and activation
    functions. The parameters \textit{in italics} vary from their original papers and open-source code. Note that a single FNO mode actually consists of \textbf{two} trigonometric
    polynomials and is indicated by ``($\mathbb C$)''.}
  \label{tab:models-para}
\begin{center}
  \begin{tabular}{c|cccccccc}
    \hline
    Hyperparameter & bs & epochs & optimizer & modes & width & $L$ & grids & $\sigma$\\ \hline
    FNO & $20$ & $\it{5000}$ & ADAM & $\mathit{20}(\mathbb C)$ & $\mathit{50}$ & $4$ & Uniform & GeLU\\ \hline
    OPNO & $20$ & $5000$ & ADAM & $40$ & $50$ &$4$ & CGL & GeLU \\ \hline
    IAE & $50$ & $\it{5000}$ & ADAM  & $256$ & $64$&$4$ &Uniform & ReLU\\ \hline
    GT & $4$ & $\it{5000}$ & 1cycle & $\mathit{20}(\mathbb C)$ & $\mathit{50}$& $4$ & Uniform & ReLU\\ \hline
    POD-DO & $20$ & $500000$ & ADAM & $32$ & -- &-- & Uniform & $\tanh$\\ \hline
  \end{tabular}
\end{center}
\end{table}
\subsection{Experiment 1: Viscous Burgers equation with Neumann BCs}\label{exp:burgers}
We consider the one-dimensional viscous Burgers equation
\begin{equation}
  \partial_t u(x, t) + \frac{1}{2}\partial_x (u^2(x,t)) = \nu \partial_{xx}u(x, t), \ x \in I \label{eq:burgers-neumann}
\end{equation}
subject to the initial-boundary conditions
\begin{eqnarray}
  &u(x, 0) = u_0(x), \\
  &\frac{\partial u}{\partial x} (\pm 1, t) = 0, \ t>0, \label{eq:bc-neumann}
\end{eqnarray}
and we aim to learn the solution operator $S(1) : S(1)u_0 \mapsto u(\cdot, 1), u_0 \in H^s(I)$. Burgers' equation is an important PDE in various fields, such as fluid dynamics, traffic flow and shock wave
theory, while the Neumann BC
represents the free flow of fluid in/out of the boundary and is quite common in its modeling scenarios. The initial condition $u_0(x)$ is generated using a Gaussian random field according to $u_0 \sim \mu$, where $\mu = \mathcal N(0,
625(-4\Delta + 25 I)^{-2})$ with Neumann BCs, its Karhunen--Lo\'{e}ve (K--L) expansion being of the form
\begin{equation}
\label{eq:39}
u_0(x) = \sum\limits_{k=0}^{\infty} 25( 4(k\pi/2)^2 + 25)^{-1} \cdot \mathrm{Re}\left\{ e^{\mathrm{i}k\pi x/2} \right\} \cdot \xi_k , \ x \in \left[ -1, 1 \right],
\end{equation}
where $\left\{ \xi \right\}_k$ are i.i.d. standard Gaussian random variables. In addition, the viscosity is set as $\nu = 0.1/
\pi$; then, $1000$ instances are generated for training as well as $100$ instances for test. Therefore, except for the effect of different BCs, such a dataset has almost the same pattern as that of the FNO dataset for Burgers' equation with
the periodic BC in \cite{li2021fourier}, and we can assert that such a comparison is fair.

The global relative $L^2$
error and the $L^{\infty}$ error of the BCs are shown in   \cref{tab:exp-1.1}, where the $L^{\infty}$ norm BC error for CGL grid outputs is computed by the polynomial
differentiation method \cref{eq:poly-diff}, and those for uniform grid outputs are approximated by the first-order difference; i.e., let $h = \frac{2}{N}$ denotes the spatial step size, then
\begin{equation}
  \mathcal E_{b.c. L^{\infty}} = \mathrm{max} \left\{ \left| \delta_x u_0 \right| , \left| \delta_x u_{N-1} \right| \right\},\nonumber
\end{equation}
where
$$\delta_x u_n \triangleq \frac{u_{n+1} - u_{n}}{h}$$

\begin{figure}[htbp]
\centering
    \subfloat[][test error in Experiment 1]{\includegraphics[width=.33\linewidth]{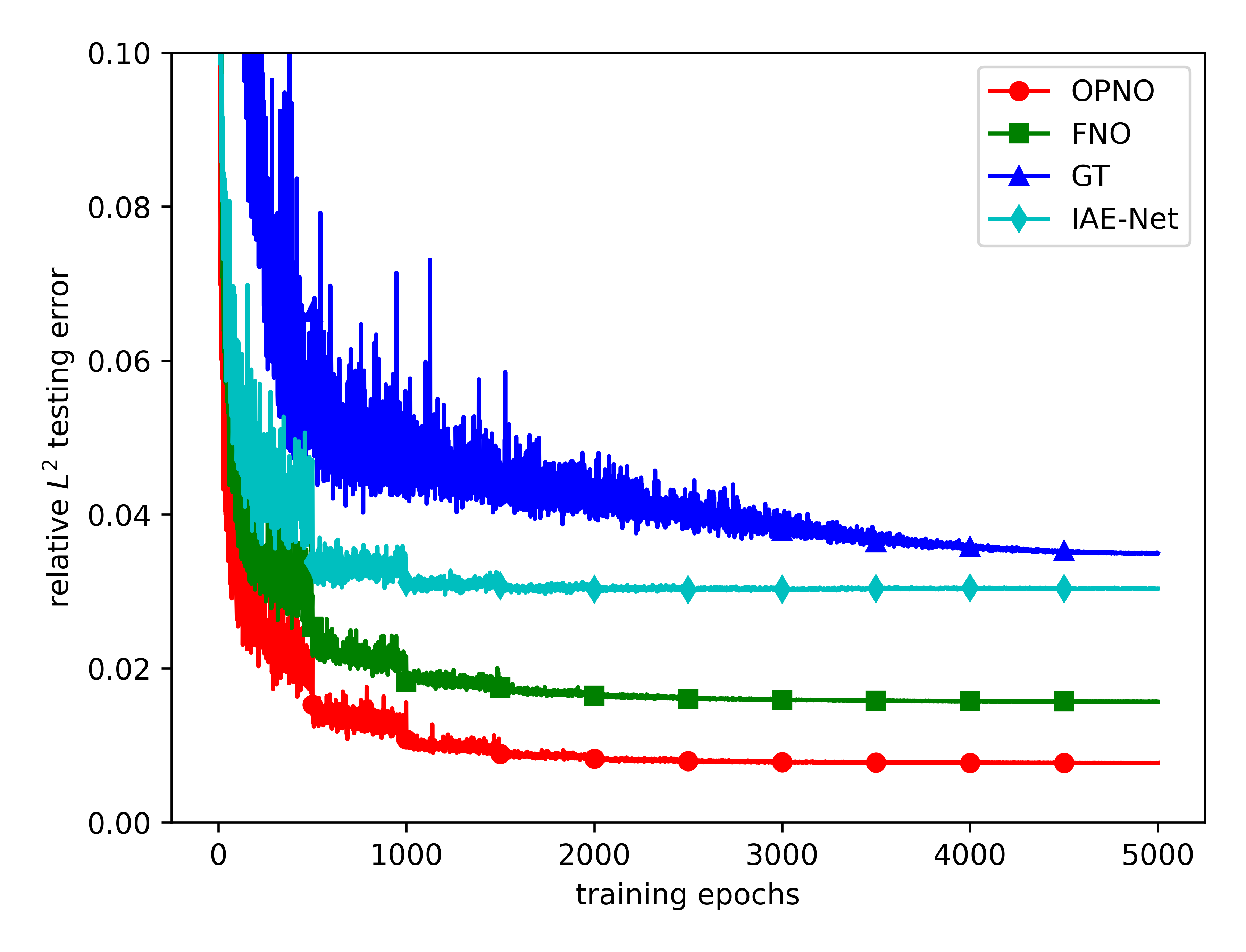}\label{fig:non-overfit}}
    \subfloat[][test error in Experiment 2]{\includegraphics[width=.33\linewidth]{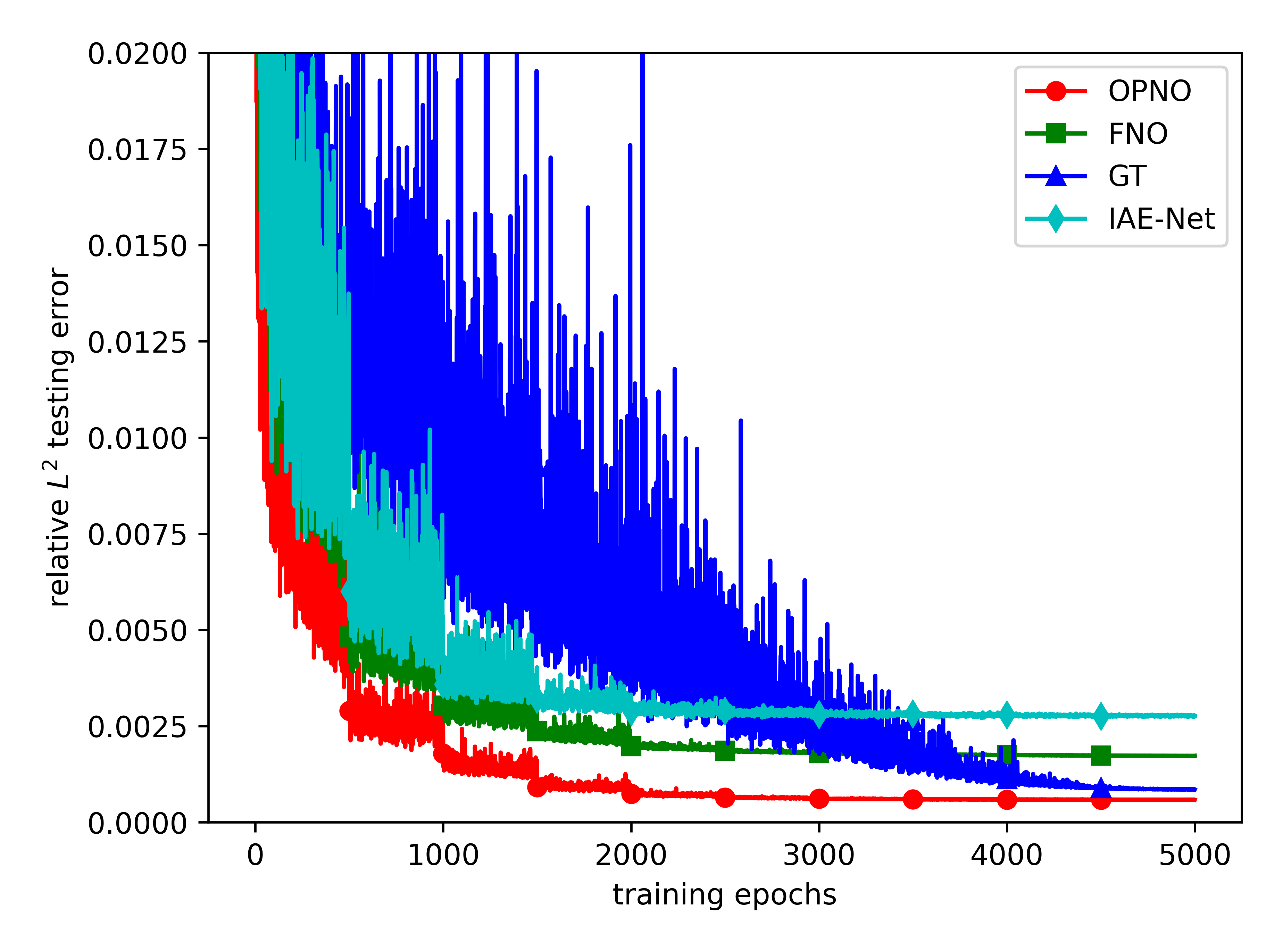}\label{fig:non-overfit2}}
    \subfloat[][quasi-linear computational complexity]{\includegraphics[width=.33\linewidth]{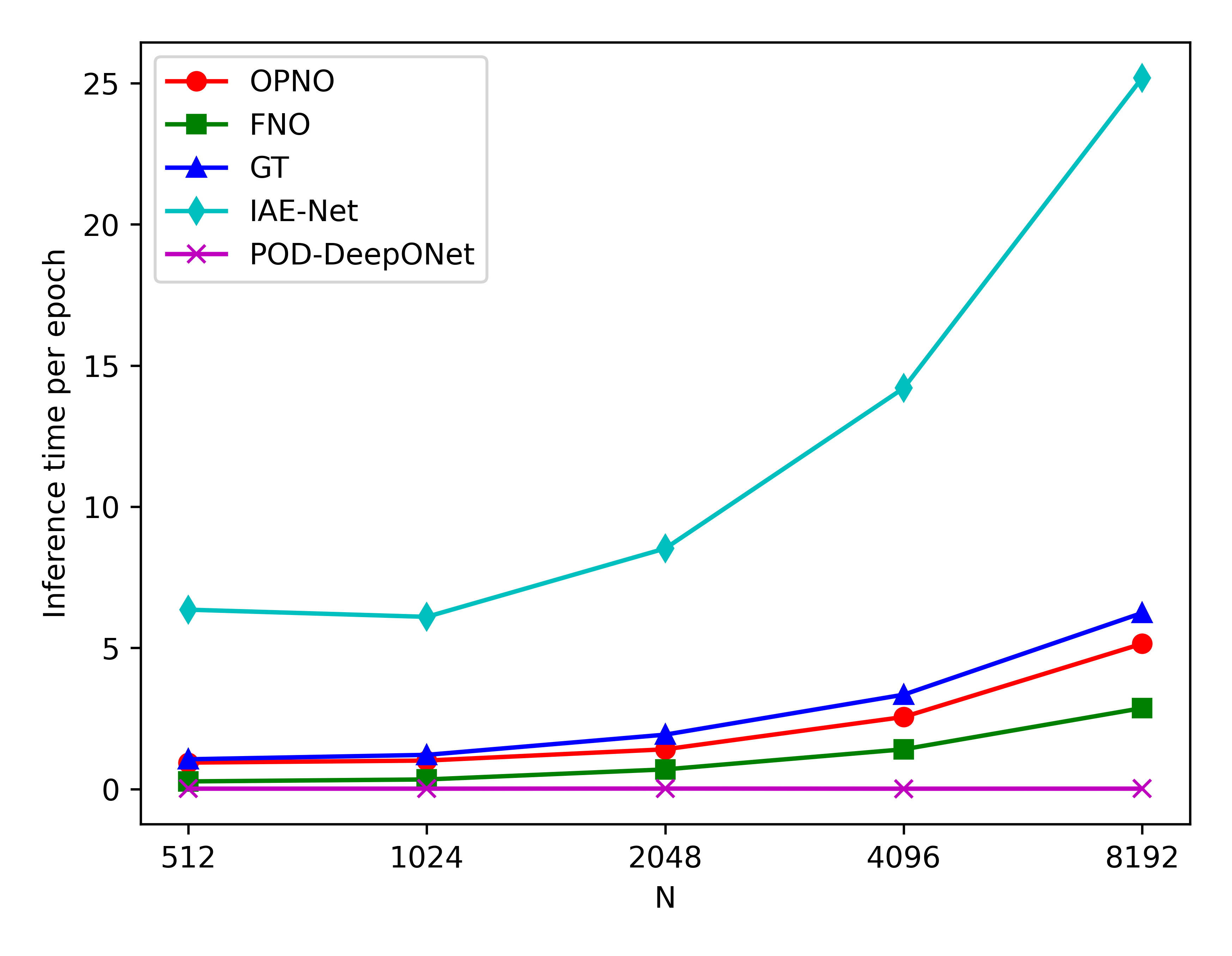}\label{fig:GPU-time}}
    \caption{\textbf{The relative $L^2$ errors with respect to the number of training iterations} in Burgers-Neumann-BC experiment (a) and heat-Robin-BC experiment (b) with $N=256$, where the neural operators exhibit non-overfitting capacities after long
training periods. \textbf{The GPU
time required per epoch} (c) in heat-Robin-BC experiment, where the batch sizes of all models are set as $20$.}
\end{figure}

\begin{figure}[htbp]
  \centerline{\includegraphics[width=0.8\textwidth]{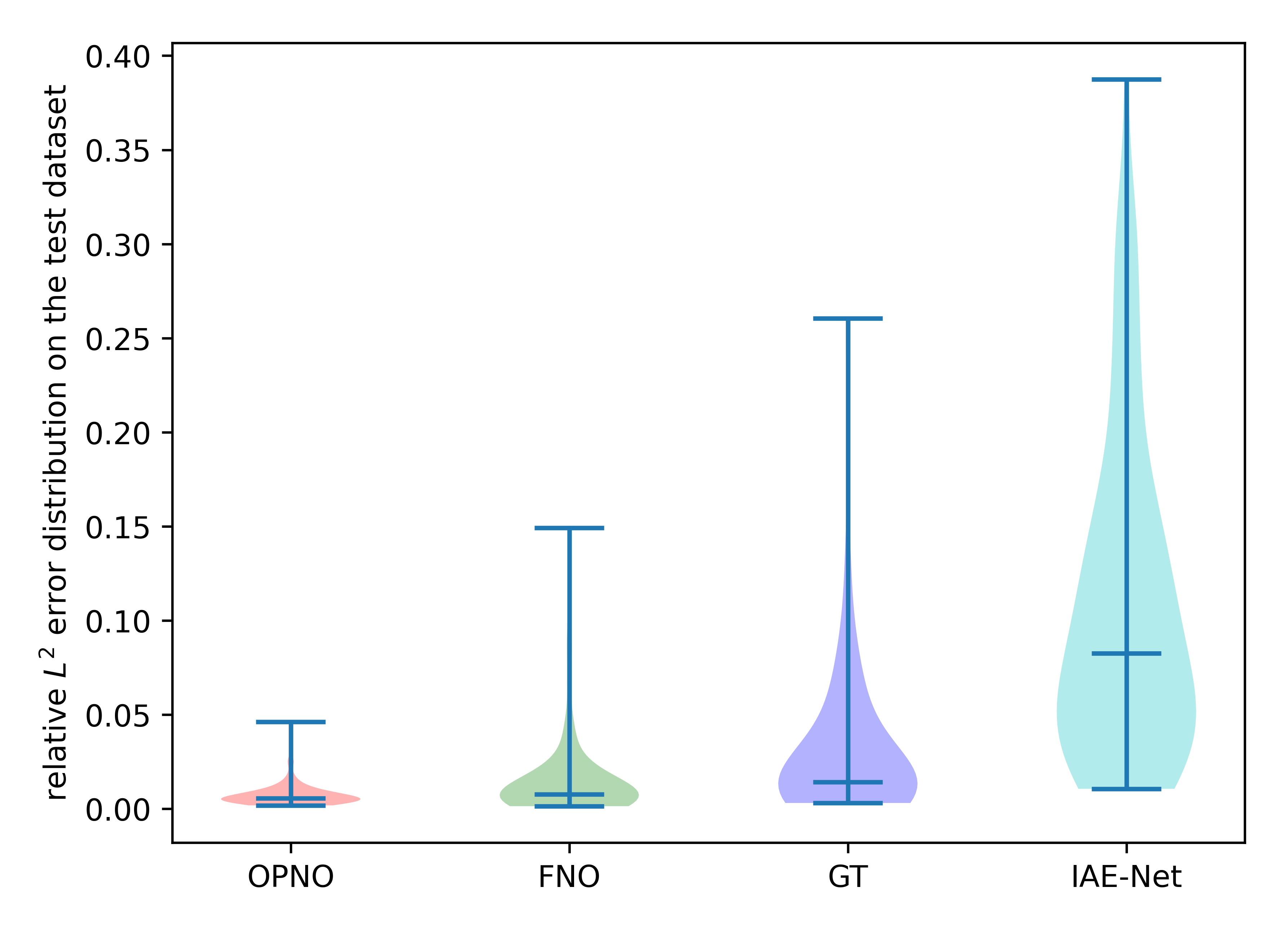}}
  \caption{ $L^2$ error distribution of the Burgers-Neumann-BC experiment.}
  \label{fig:violin}
\end{figure}

From the results, we can deduce the following conclusion.
\begin{itemize}
\item As   \cref{tab:exp-1.1} shows, The OPNO outperforms other models not only in terms of the BC accuracy but also in terms of the global precision. Further evidence is shown in  \cref{fig:neumann} that the OPNO generate
  predictions that are ``globally accurate''.
\item It is shown in  \cref{fig:violin} that the higher average errors of non-BC-satisfying models result from the bad performance on a handful of ``hard'' samples, which implies
the OPNO is more reliable in terms of the worst-case performance due to the BC satisfaction property.
\item  \cref{fig:non-overfit} illustrates that the test errors of FNO, OPNO, GT, and IAE-Net do not increase after long-term training, so these models possess the self-regulation properties when solving the PDEs we are interested in, and the decision of extending the number of training epochs successfully improved
their performance by mitigating underfitting. Other experiments \textbf{also} follow the same non-overfitting pattern, see \cref{fig:non-overfit2} for an example.
\begin{table}[htbp]
  \footnotesize
  \caption{Benchmarks on 1D Burgers equation with Neumann BCs}
  \label{tab:exp-1.1}
\begin{center}
\resizebox{1\columnwidth}{!}{
  \begin{tabular}{ccc|cc|cc}
    \hline
     N &\multicolumn{2}{c|}{256} &\multicolumn{2}{c|}{1024}&\multicolumn{2}{c}{4096}\\ \cline{2-7}
  &$L^2$&b.c. $L^{\infty}$ &$L^2$&b.c. $L^{\infty}$ &$L^2$&b.c. $L^{\infty}$ \\\hline
    FNO & $1.571e-2$ & $2.902e-1$ &$1.684e-2$&$4.092e-1$  &$1.688e-2$&$5.478e-1$\\
    \hline
    OPNO & $\mathbf{7.704e-3}$&$\mathbf{5.963e-12}$ &$\mathbf{7.814e-3}$&$\mathbf{1.131e-10}$ & $\mathbf{7.821e-3}$&$\mathbf{1.900e-9}$\\ \hline
    IAE & $3.285e-2$ & $1.056$ & $2.842e-2$ & $3.494$ & $2.660e-2$ & $2.016e+1$ \\ \hline 
    POD-DO & $ 1.210e-1$ & $ 1.285e-1$ & $1.337e-1$ & $4.731e-2$ & $1.410e-1$& $1.093e-2$ \\ \hline
    GT & $3.744e-2$ & $1.187$ & $3.902e-2$ &$1.411$  &$3.848e-2$ & $1.791$\\ \hline
    %OPNO(width$=150$) & $1.761e-2}$ & $2.440e-1}$ & $1.747e-2}$ &$2.780e-1}$  &$1.781e-2}$ & $3.616e-1}$\\ \hline
  \end{tabular}
}
\end{center}
\end{table}

%\item Compared with the FNO, although the OPNO is more memory-costly regarding the frequency interaction when the bandwidth $w>1$, the additional parameter demand of a constant factor is worthwhile since devoting resources to the frequency domain (for example, increasing the width of frequency features of FNO) can \textit{not} compensate for the gap in accuracy. See the FNO (with$=150$) result in   \cref{tab:exp-1.1}.
\end{itemize}
\subsection{Experiment 2: Heat diffusion equation with Robin BCs}\label{exp:robin}
When referring to the heat diffusion equation, the Robin BC is also known as the convection BC, which describes the convection heating/cooling occurring at an object's surface and is crucial in many engineering fields,
including the designs of chips and engines. Consider the equations
\begin{eqnarray}
%\label{eq:54}
  &\partial_t u(x, t) - k u_{xx} + F(u, x) = 0, \ x \in I, \label{eq:38} \\
  &u(x, 0) = u_0(x), \ x \in I, \\
  &C_1 u(-1, t) - k_T u_x(-1, t) = C_2, \ C_3 u(1, t) + k u_x(1, t) = C_4, \ t \geq 0, \label{eq:robin-bc}
\end{eqnarray}
where $k_T$ represents the conduction heat flux, and we are to investigate the approximation of $\mathcal G: \mathcal G(u_0(x)) = u(x, 1)$.

Furthermore, we remain curious about whether and how the test errors of neural operators can drop to significantly low levels in practice, but doubtlessly, a considerably large training dataset is
helpful. To generate such a dataset (up to $10^6$ training samples and $1000$ test samples) using numerical methods, we fix the equation parameters $k_T=C_1=C_3=0.02$,
$C_2=C_4=0$ and $F(u, x) \equiv 0$. To satisfy the BCs and provide the dataset with sufficient degrees of freedom, analogously to Eq. \cref{eq:39}, the initial condition $u_0(x)$ is generated by a quasi-polynomial chaos expansion:
\begin{equation}
\label{eq:40}
u_0(x) = \sum\limits_{k=0}^{\infty}  \hat v_k\phi_k(x) \xi_k, \nonumber
\end{equation}
where $\phi_k(x) = T_k(x) - \frac{k-1}{k+3} T_{k+2}(x)$ satisfies the Robin BC \eqref{eq:robin-bc}, and
$$\hat v_k = \frac{2\sigma}{\pi p_k}  \gamma^{-k}. $$
We further set $\gamma = 1.1$, $\sigma=4$, and $p_k = \tilde c_k + \left( \frac{k^2+1}{(k+2)^2+1} \right)^2$.

\begin{figure}[htbp]
\centering
    \subfloat[][Example \uppercase\expandafter{\romannumeral1}]{\includegraphics[width=.33\linewidth]{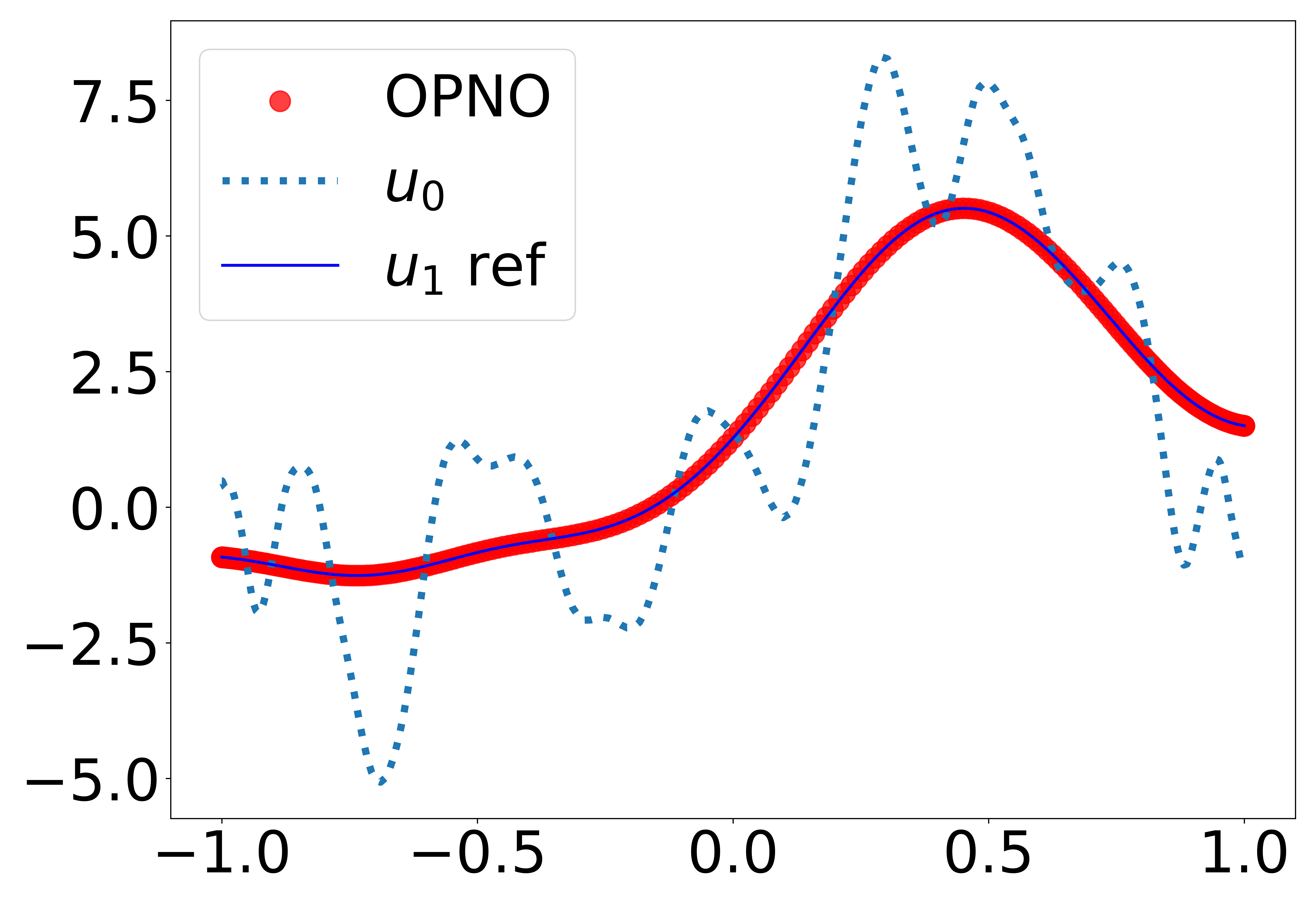}}
    \subfloat[][Example \uppercase\expandafter{\romannumeral2}]{\includegraphics[width=.33\linewidth]{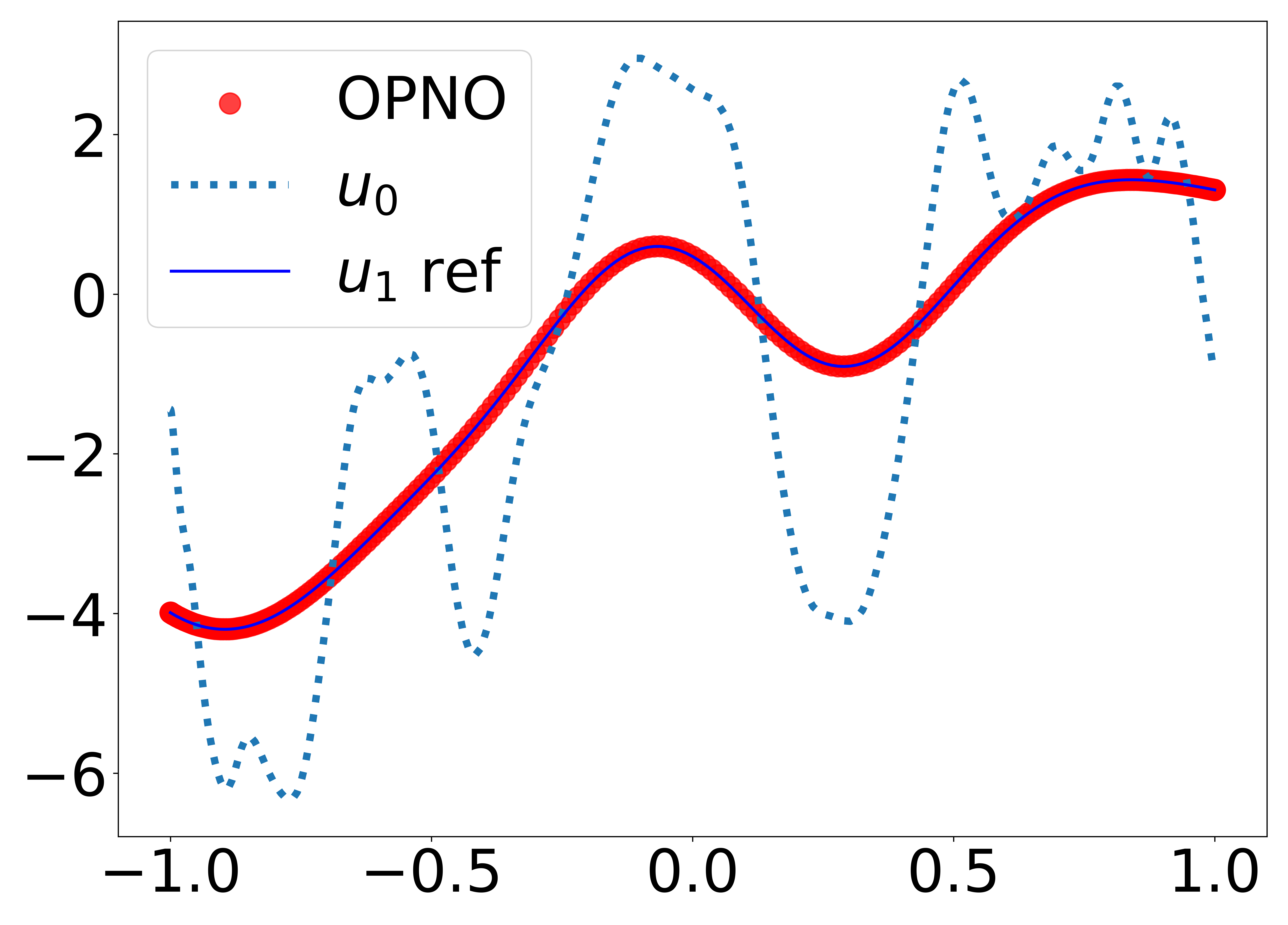}}
    \subfloat[][Example \uppercase\expandafter{\romannumeral3}]{\includegraphics[width=.33\linewidth]{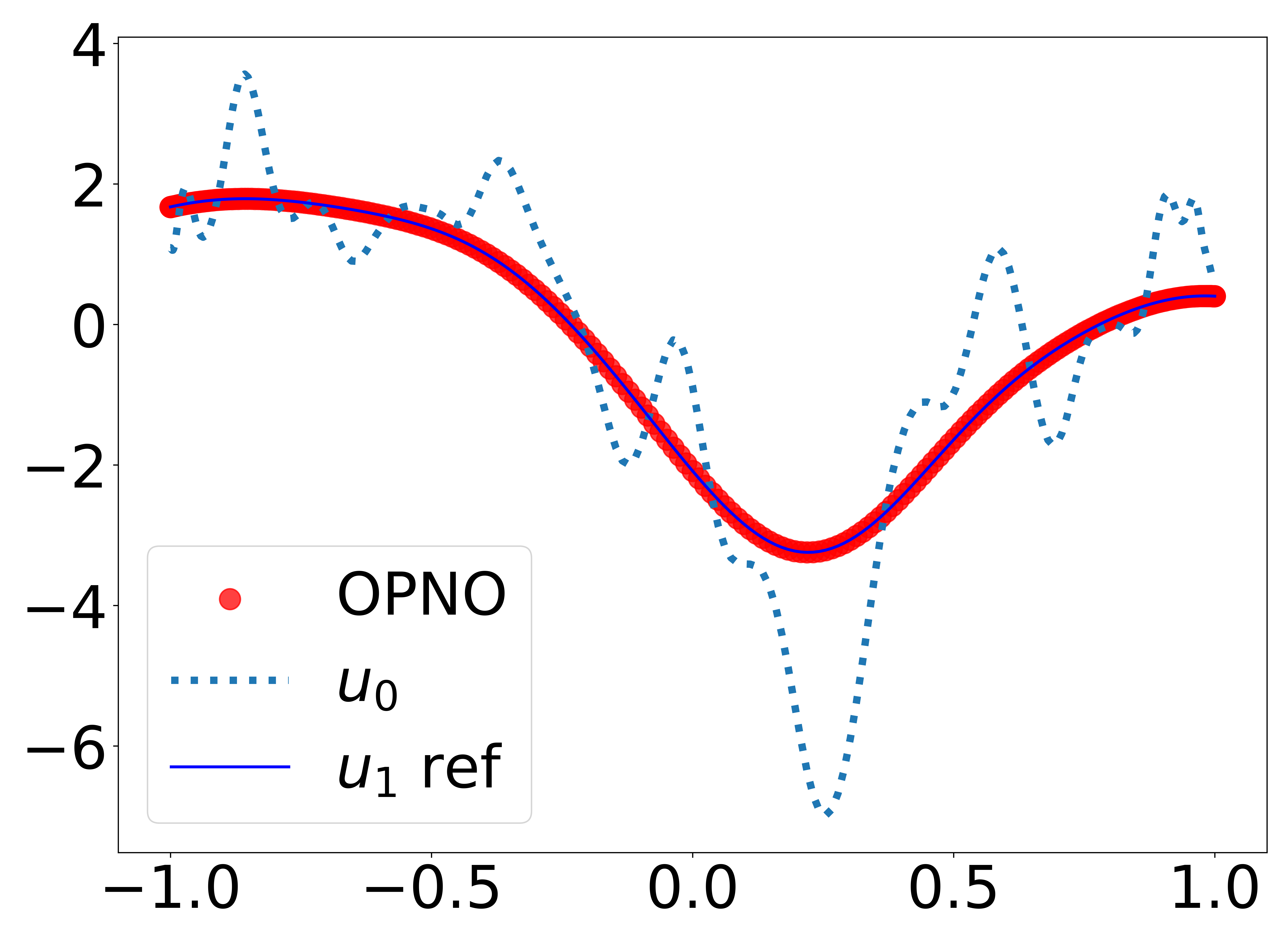}}
    \caption{Evaluation results obtained by the OPNO for 3 test instances of the heat-Robin-BC experiment.}
   \label{fig:robin}
\end{figure}

\textbf{Comparison between deep learning methods}: First, we perform an experiment on the sub-training set consisting of $1000$ training instances, and the results are shown in
  \cref{tab:exp-2.1}. Although the native OPNO seems to perform slightly worse than the GT, we note that the mini-batch strategy is important for neural operator training, and the OPNO overtakes the GT model
once both models adopt the same batch size. Therefore, taking all factors into account, the OPNO obtains the lowest relative
error among the tested models, and the output accurately satisfies the Robin BCs within machine precision. Preditions for three test instances are demonstrated in  \cref{fig:robin}, where the
prediction and reference solutions coincide.

We also test the GPU time usage (seconds per epoch) required for training different models, and we can determine from  \cref{fig:GPU-time} that the OPNO model runs in
quasi-linear time, whereas the POD-DeepONet is the fastest tested neural operator.
\begin{table}[htbp]
  \footnotesize
  \caption{Benchmarks on 1D heat equation with Robin BCs.}
  \label{tab:exp-2.1}
  \begin{center}
\resizebox{1\columnwidth}{!}{
  \begin{tabular}{ccc|cc|cc}
    \hline
     N &\multicolumn{2}{c|}{256} &\multicolumn{2}{c|}{1024}&\multicolumn{2}{c}{4096}\\ \cline{2-7}
  &$L^2$&b.c. $L^{\infty}$ &$L^2$&b.c. $L^{\infty}$ &$L^2$&b.c. $L^{\infty}$ \\\hline
    FNO &$1.880e-3$&$6.030e-1$ & $1.728e-3$ & $8.394e-1$ & $1.722e-3$ &$1.012$\\
    \hline
    OPNO &$6.073e-4$&$2.325e-11$ & $6.585e-4$&$4.595e-10$ & $6.751e-4$&$7.055e-9$\\ \hline
    IAE& $2.333e-3$ & $1.307$ & $1.618e-3$ & $2.485$ & $1.610e-3$ & $6.751$\\ \hline
    POD-DO& $4.052e-1$ & $1.832e-1$ & $4.156e-1$ & $4.949e-2$ & $4.084e-1 $ & $1.956$ \\ \hline
    \makecell{GT\\($bs=4$)} & $9.400e-4$&$5.634e-1$ & $6.058e-4$ & $9.203e-1$ & $6.830e-4$ & $1.491$ \\ \hline
    \makecell{OPNO\\($bs=4$)} & $\mathbf{2.927e-4}$&$\mathbf{2.178e-11}$ & $\mathbf{2.423e-4}$&$\mathbf{4.234e-10}$ &$\mathbf{2.843e-4}$ & $\mathbf{6.888e-9}$\\ \hline
  \end{tabular}

}
\end{center}
\end{table}

\textbf{Comparison with the numerical method}: Next, we compare OPNO with the traditional numerical method in terms of not only computational efficiency but also accuracy. Numerical methods with 2nd-order convergence are the most popular PDE solvers in
scientific software due to their balance between efficiency and stability. To be precise, for Eq. \eqref{eq:38}, we adopt the centered second difference in the spatial direction and the second-order
Crank--Nicolson method in the temporal direction, while the ghost
point method is applied on the boundary to discretize the BCs \eqref{eq:robin-bc}. Let $h$ and $\tau$ denote the spatial and temporal step sizes, respectively, and fix $h = \tau$.
%The 2nd-order FDM is of the form
We recompute the $1000$ test instances using
this 2nd-order FDM, and the errors and the average time consumption levels obtained with different $h, \ \tau$ are listed in   \cref{tab:fdm}. The results verify the 2nd-order convergence rate of FDM. 

Note that system \eqref{eq:38}--\eqref{eq:robin-bc} is already one of the most trivial problems for the FDM to compute since only a sparse system of linear equations needs to be solved at each time
step. As the mesh is refined, the discretization errors of the FDM decrease to an arbitrarily small value, but the
computational cost also significantly increases. Nevertheless, the inference time of the OPNO is independent of the size of the training set, whereas the self-regularization feature suggests that the test
errors of neural operators may reach desirably low levels once trained on a sufficiently large training set.

\begin{figure}[tbhp]
\centerline{\includegraphics[width=1\textwidth]{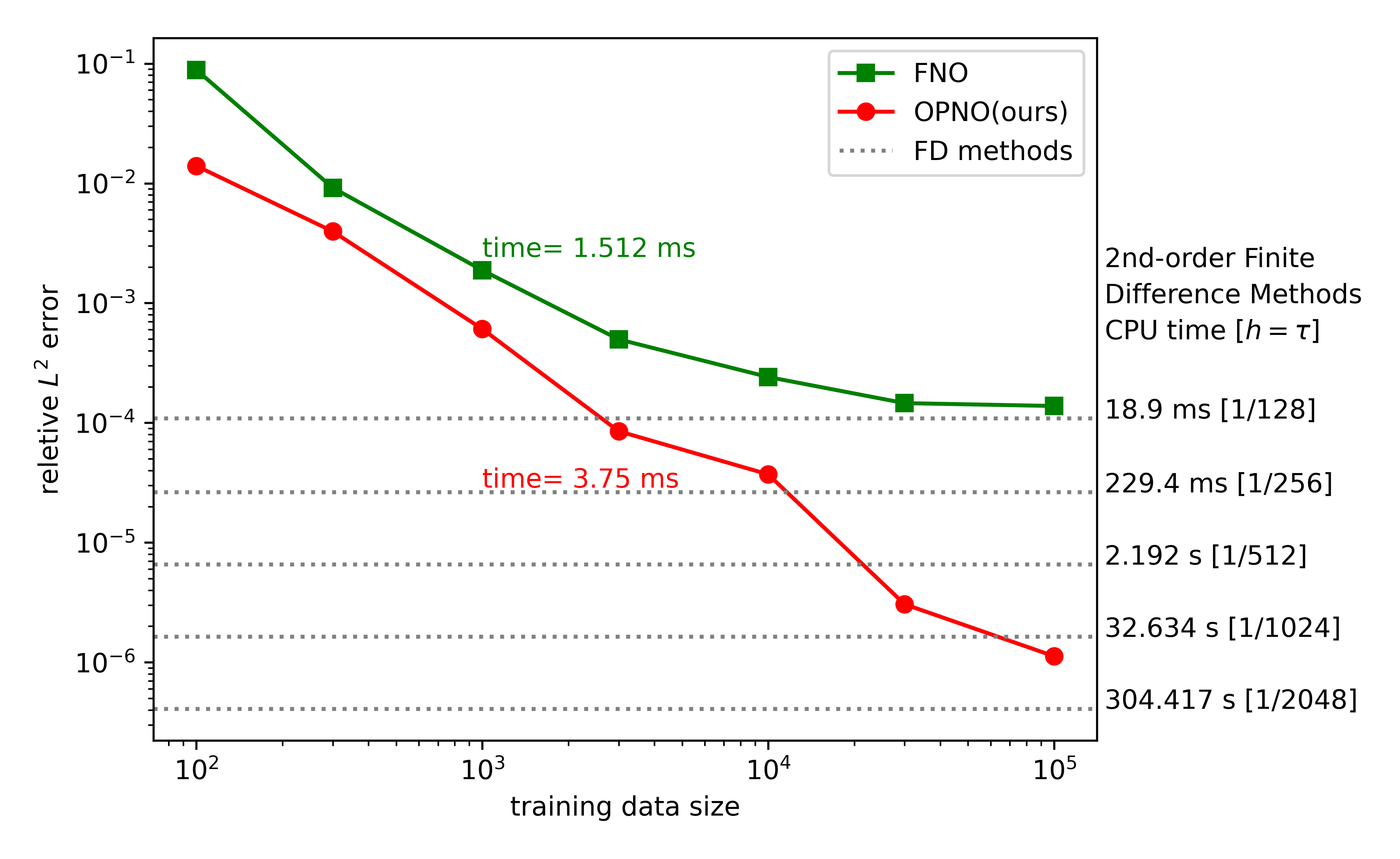}\label{fig:vs-fdm}}
\caption{The inference times and relative $L^2$ errors of neural operators versus those of the 2nd-order FDM. For the neural operators, $N=256$, and the batch size is set as $20$.}
\end{figure}

Consequently, we test the FNO and OPNO on sub-training sets with sizes ranging from
$10^{2}$ to $10^{6}$, and the relative $L^2$ errors are listed in   \cref{tab:exp2-datasize}. A visual comparison chart against 2nd-order FDM is shown in  \cref{fig:vs-fdm}. The OPNO model has an inference time of only
$3.75$ ms compared to the $304,417$ ms average computational time of the FDM and achieves a lower error than the 2nd-order
FDM with a reasonably fine mesh possessing $h, \tau< 0.001$. While the results are already impressive, we emphasize that the error may be further reduced if better training techniques are employed. This also demonstrates the importance of
choosing the correct underlying SOL basis that satisfies the given BCs.

\begin{table}[htbp]
  \footnotesize
\begin{center}
\caption{The average relative $L^2$ errors, numerical convergence order and time consumption (seconds) per instance of the \textbf{2nd-order FDM}.}\label{tab:fdm}
\begin{tabular}{c|ccccc}
  \hline
  $h, \tau$ & $1/128$ & $1/256$ & $1/512$ & $1/1024$ &$1/2048$ \\ \hline
  $L^2$ error & $1.091e-4$ &$2.635e-5$ &$6.563e-6$ &$\mathbf{1.640e-6}$ &$4.101e-7$ \\ \hline
  order & -- &$2.050$ &$2.005$ &$2.001$ &$2.0000$ \\ \hline
time (sec.) & $0.0189$ & $0.2294$ &$2.192$ &$32.634$ &$\mathbf{304.417}$ \\ \hline 
\end{tabular}
\end{center}
\end{table}
\begin{table}[htbp]
  \footnotesize
  \caption{Relative $L^2$ errors induced by \textbf{neural operators} on training datasets with different sizes. $N=256$, bs$=20$.}
  \label{tab:exp2-datasize}
\begin{center}
    \begin{threeparttable}
  \begin{tabular}{c|cccc}
    \hline
     training data& 100 & 300 &1000 & 3000\\ \hline
    FNO & $8.858e-2$ & $9.175e-3$&$1.880e-3$ &$4.978e-4$\\ \hline
    OPNO & $1.397e-2$ & $3.964e-3$ & $6.073e-4$ & $8.509e-5$\\ \midrule[1pt]
    % OPNO-cos & & & $5.058e-4}$ &5.953e-5} & & \\ \hline
    % OPNO(1cycle) & & & $3.335e-4}$ &6.016\times 10^-5 & & \\ \hline
     training data &10000 & 30000 & 100000\\ \hline
    FNO & $2.412e-4$ & $1.460e-4$ & $1.382e-4$\\ \hline
    OPNO & $3.711e-5$ & $3.054e-6$(*) & $1.124e-6$(*)\\ \hline
  \end{tabular}
 \begin{tablenotes}
        \footnotesize
        \item[(*)] As the data size increases, the weight decay is set as $0$, and each model is trained for $6000$ epochs to mitigate underfitting.
      \end{tablenotes}
  \end{threeparttable}
\end{center}
% 2.9373e-5
\end{table}

\subsection{Experiment 3: Heat diffusion equation with inhomogeneous Dirichlet BCs}
Although solving PDEs with Dirichlet BCs is usually quite straightforward, this experiment aims to introduce a specific approach for solving PDEs with inhomogeneous BCs using the OPNO. Consider the
approximation for the solution operator $S(1) : S(1) u_0 = u(x, 1)$ of the following
heat diffusion equation with Dirichlet BCs:
\begin{eqnarray}
%\label{eq:56}
 & \partial_t u(x, t) - k_T u_{xx}= 0, \ x \in I,  \nonumber\\
 & u(x, 0) = u_0(x), \ x \in I, \nonumber  \\
  &   u(-1, t) = a, \ u(1, t) = b, t > 0, \label{eq:49}
\end{eqnarray}
with $k_T$ set as $0.02$. The above equations describe the heat conduction under the condition of constant surface temperature. The initial condition $u_0(x) = \Psi_0(x) + v_0(x)$ is generated according to $v_0 \sim \mu$, where $\mu = \mathcal N(0, 625(-4\Delta + 25 I)^{-2})$ with homogeneous Dirichlet BCs, and
$$ \Psi_0(x) = \frac{b-a}{2} x + \frac{a+b}{2}.$$
The BC parameters are fixed as $a = 0.3, \ b = -0.5$, and the dataset consists of $1000$ training instances and $100$ test instances.

Regarding the OPNO, the approach for satisfying the inhomogeneous BCs is to expand the basis. For the BCs \eqref{eq:49}, the polynomial basis is
expanded to $\left\{ \phi_k \right\}_{k \in \mathbb N} \cup \left\{ \Psi_0 \right\}$. The OPNO outperforms other neural operators, as shown in   \cref{tab:exp-3}. It is worth
mentioning that the POD-DeepONet model automatically learns the Dirichlet BCs accurately.
\begin{table}[htbp]
  \footnotesize
  \caption{Benchmarks on 1D heat diffusion equation with inhomogeneous Dirichlet BCs.}
  \label{tab:exp-3}
\begin{center}
\resizebox{1\columnwidth}{!}{
  \begin{tabular}{ccc|cc|cc}
    \hline
     N &\multicolumn{2}{c|}{256} &\multicolumn{2}{c|}{1024}&\multicolumn{2}{c}{4096}\\ \cline{2-7}
  &$L^2$&b.c. $L^{\infty}$ &$L^2$&b.c. $L^{\infty}$ &$L^2$&b.c. $L^{\infty}$ \\\hline
    FNO &$2.402e-4$&$2.464e-4$ & $2.364e-4$ & $4.242e-4$   &$2.360e-4$&$4.111e-4$\\
    \hline
    OPNO &$\mathbf{1.342e-4}$&$\mathbf{1.363e-15}$ & $\mathbf{1.099e-4}$ & $\mathbf{1.356e-15}$  & $\mathbf{1.115e-4}$&$\mathbf{1.362e-15}$\\ \hline
    POD-DO& $8.356e-3$ & $6.351e-12$ & $3.604e-2$ & $2.145e-13$ & $9.055e-2$ & $5.069e-12$\\ \hline
    IAE & $1.002e-3$ & $7.265e-3$ & $8.638e-4$ & $7.203e-3$ & $8.986e-4$ & $9.416e-3$ \\ \hline
    GT & $5.877e-4$ & $5.616e-4$ & $5.931e-4$ & $7.166e-4$ & $6.149e-4$ & $6.394e-4$ \\ \hline
    
  \end{tabular}
}
\end{center}
\end{table}

\begin{figure}[htbp]
\centering
    \subfloat[][Example \uppercase\expandafter{\romannumeral1}]{\includegraphics[width=.33\linewidth]{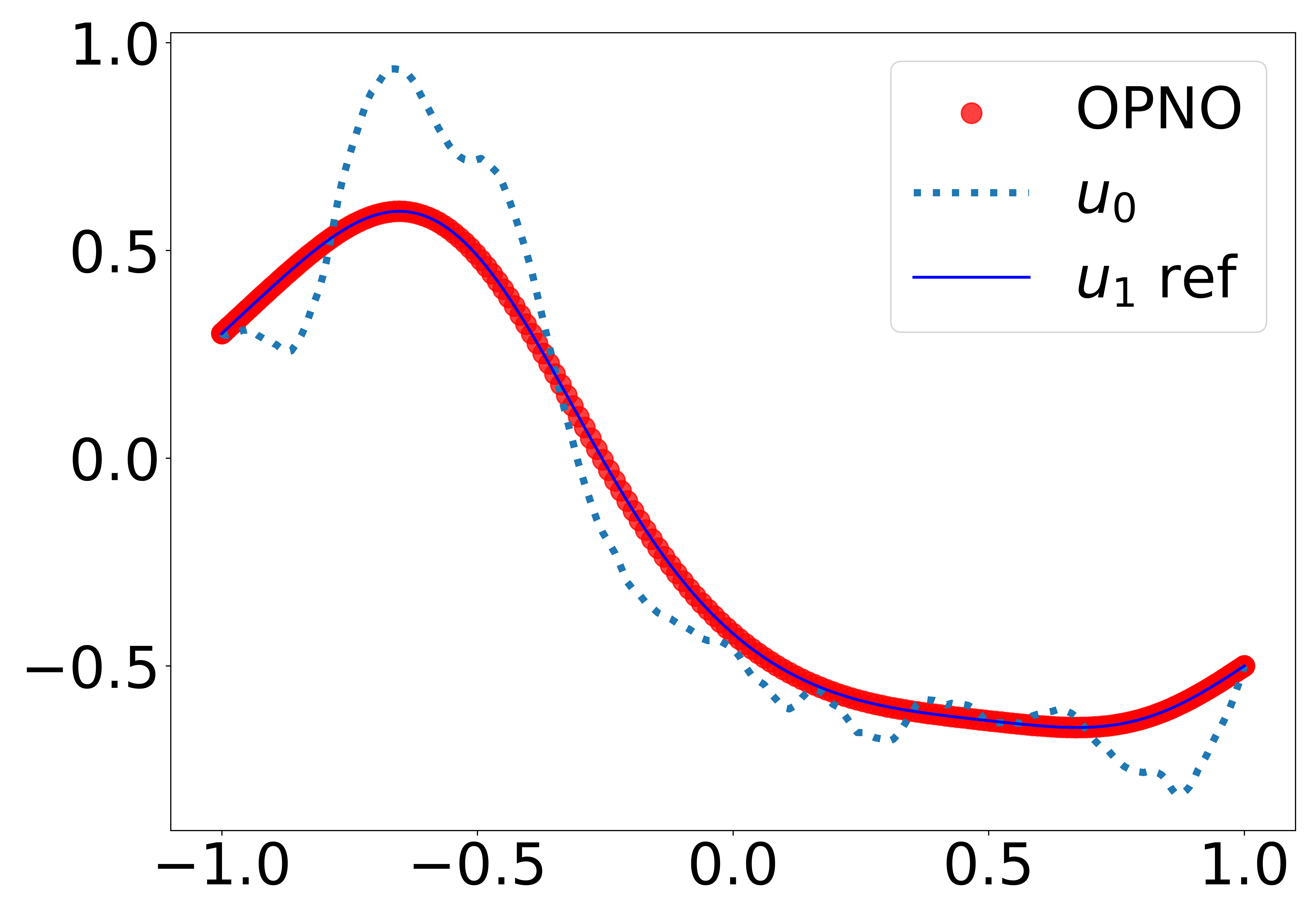}}
    \subfloat[][Example \uppercase\expandafter{\romannumeral2}]{\includegraphics[width=.33\linewidth]{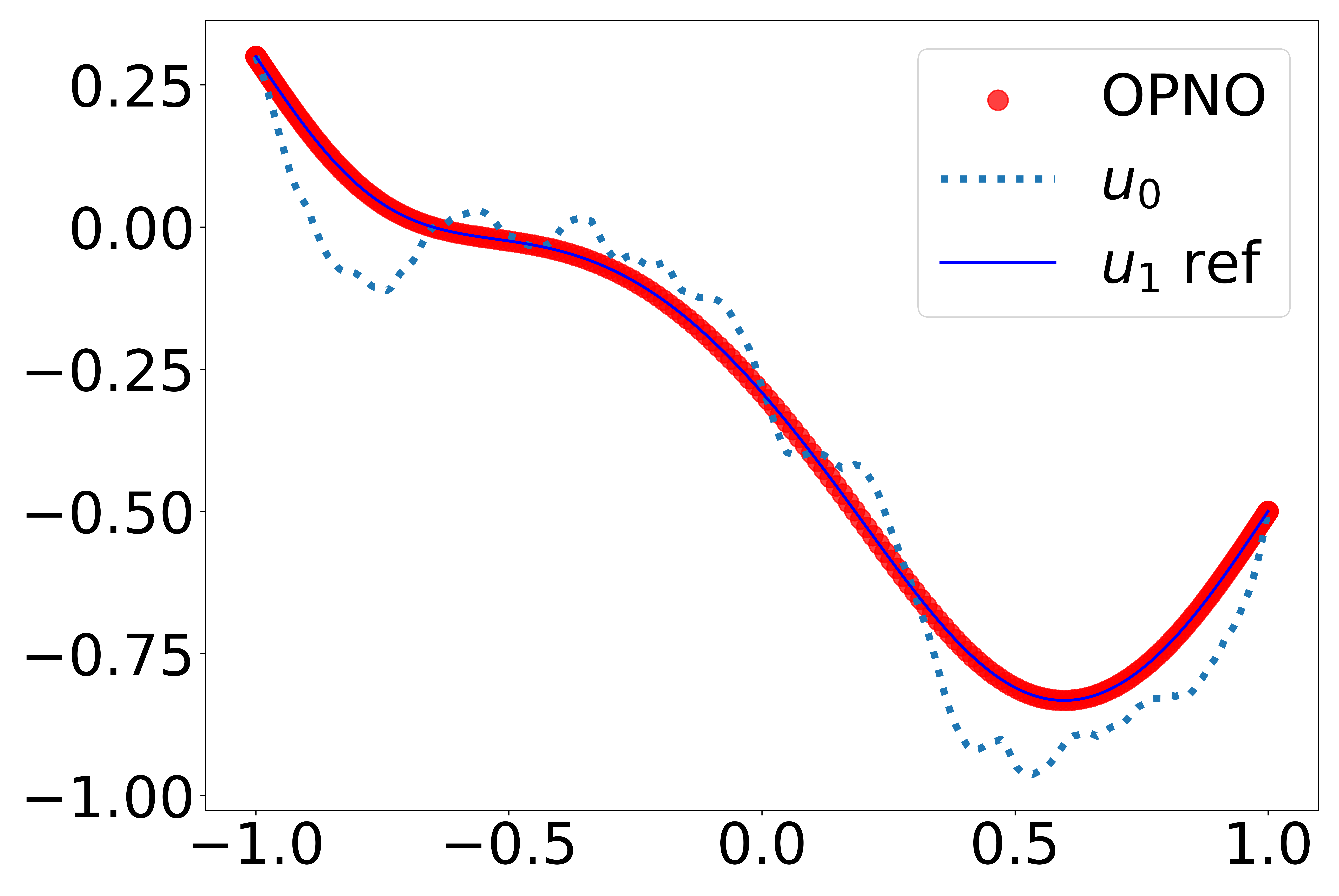}}
    \subfloat[][Example \uppercase\expandafter{\romannumeral3}]{\includegraphics[width=.33\linewidth]{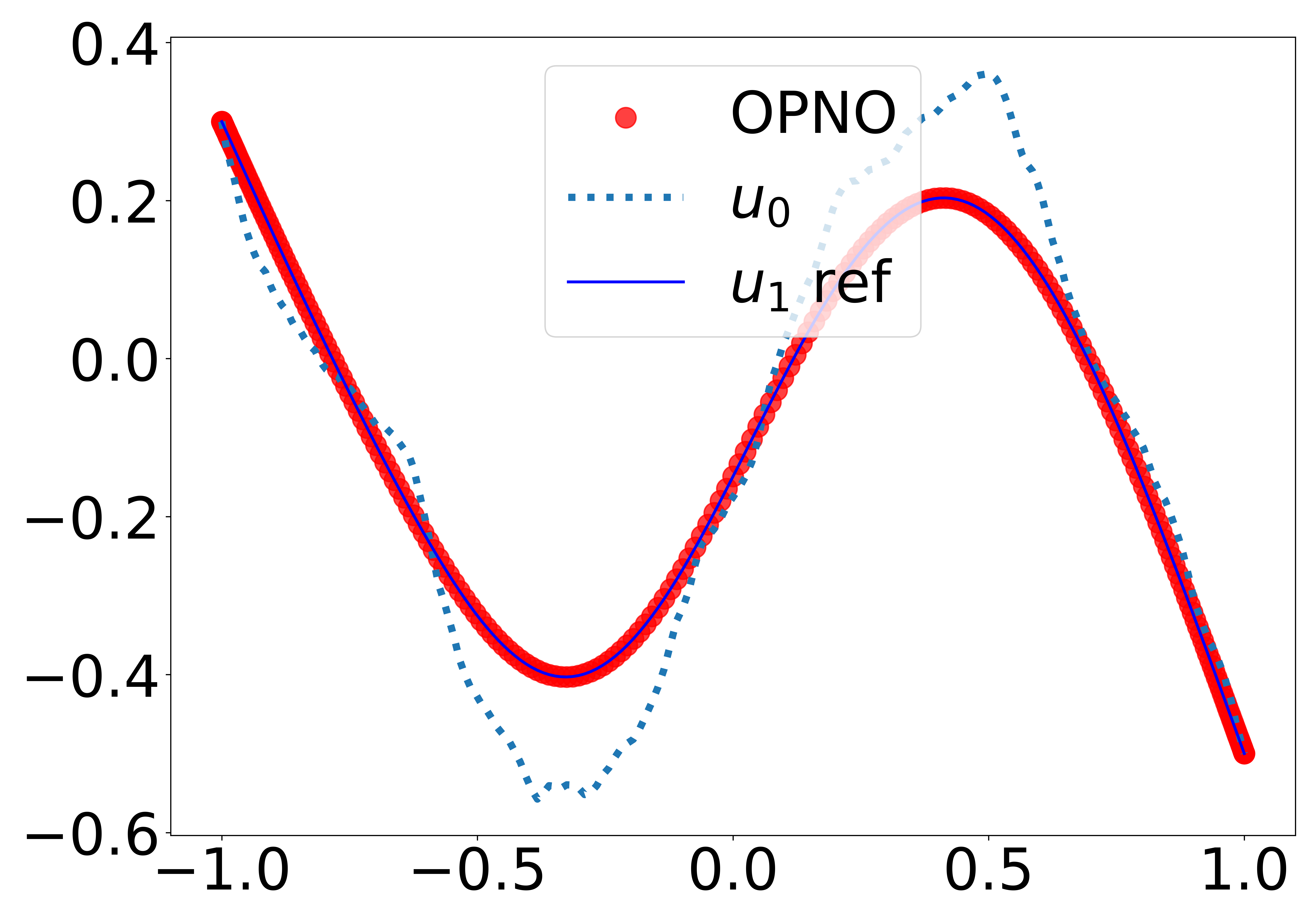}}
    \caption{Evaluation results obtained by the OPNO for 3 test instances of the heat-Dirichlet-BC experiment.}
   \label{fig:heat-predict}
\end{figure}

\subsection{Experiment 4: 2D Burgers' equation with Neumann BCs}
We consider the 2D Burgers' equation with Neumann BCs:
\begin{equation}
\label{eq:51}
\partial_t u(x, t) + (u\cdot\nabla) u = \nu \Delta u(x, t), \ x \in D \nonumber
\end{equation}
subject to
\begin{eqnarray}
  & &u(x, 0) = u_0(x), \ x \in D, \nonumber \\
  & &\frac{\partial u}{\partial \mathbf n}(x, t) = 0, \ x \in \partial D, \ t >0.\nonumber
\end{eqnarray}
where $D = I^2$. Rather than constructing a time-stepping recurrent solver, we are interested in learning the operator for directly mapping $u_0(x)$ to the set of solutions at a specific time $\left\{ u(x, t_k)
\right\}, \ t_k \in T$. We further fix $\nu = 0.001$ and $T = \left\{0.2, 0.6, 1  \right\}$. The initial conditions are generated according to $u_0 \sim \mathcal N(0, 16(\Delta + 16I)^{-2})$ with
homogeneous Neumann BCs.

\begin{figure}[tbhp]
  \begin{center}
    \subfloat[][Example \uppercase\expandafter{\romannumeral1}]{\includegraphics[width=.75\linewidth]{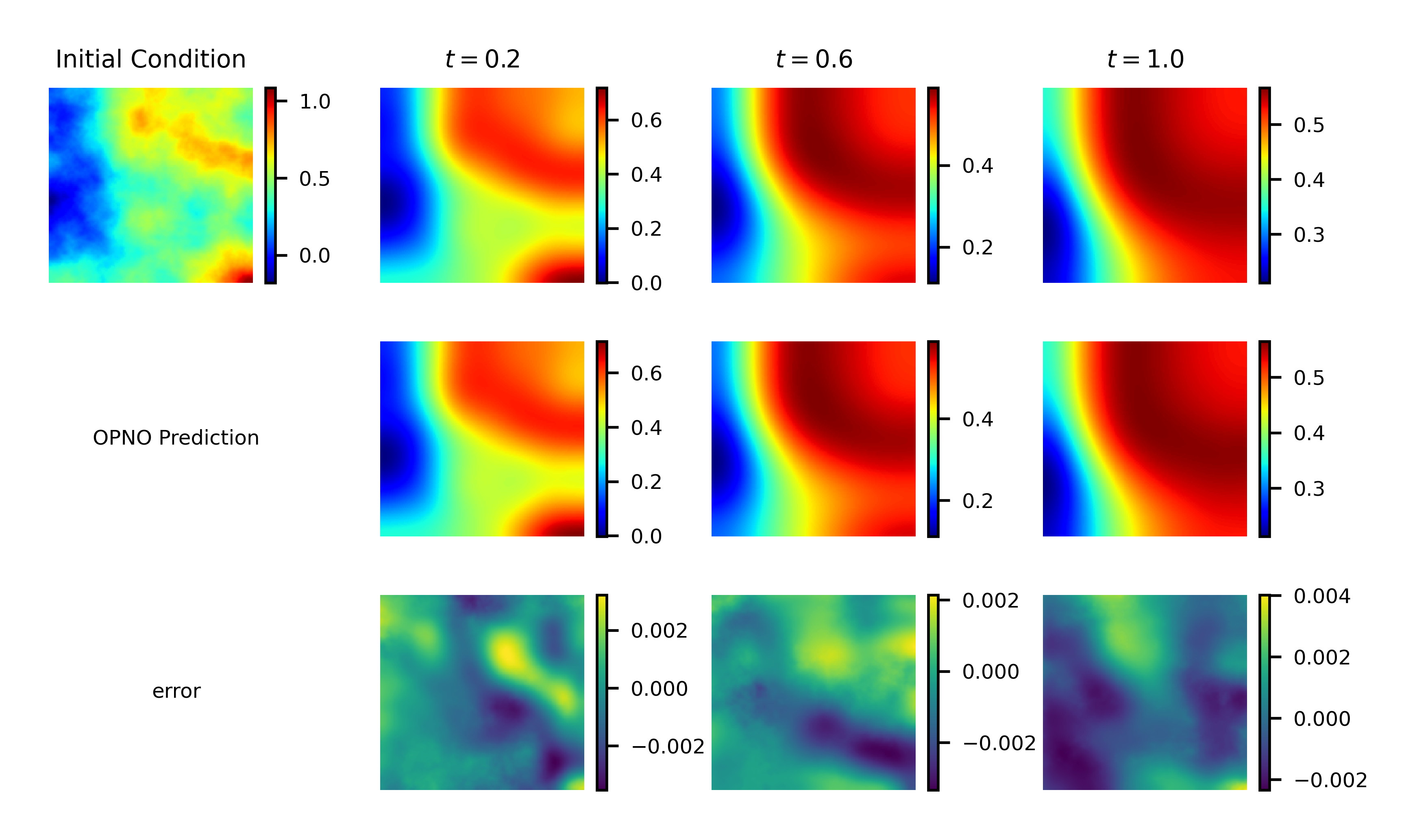}}\\
    \subfloat[][Example \uppercase\expandafter{\romannumeral2}]{\includegraphics[width=.75\linewidth]{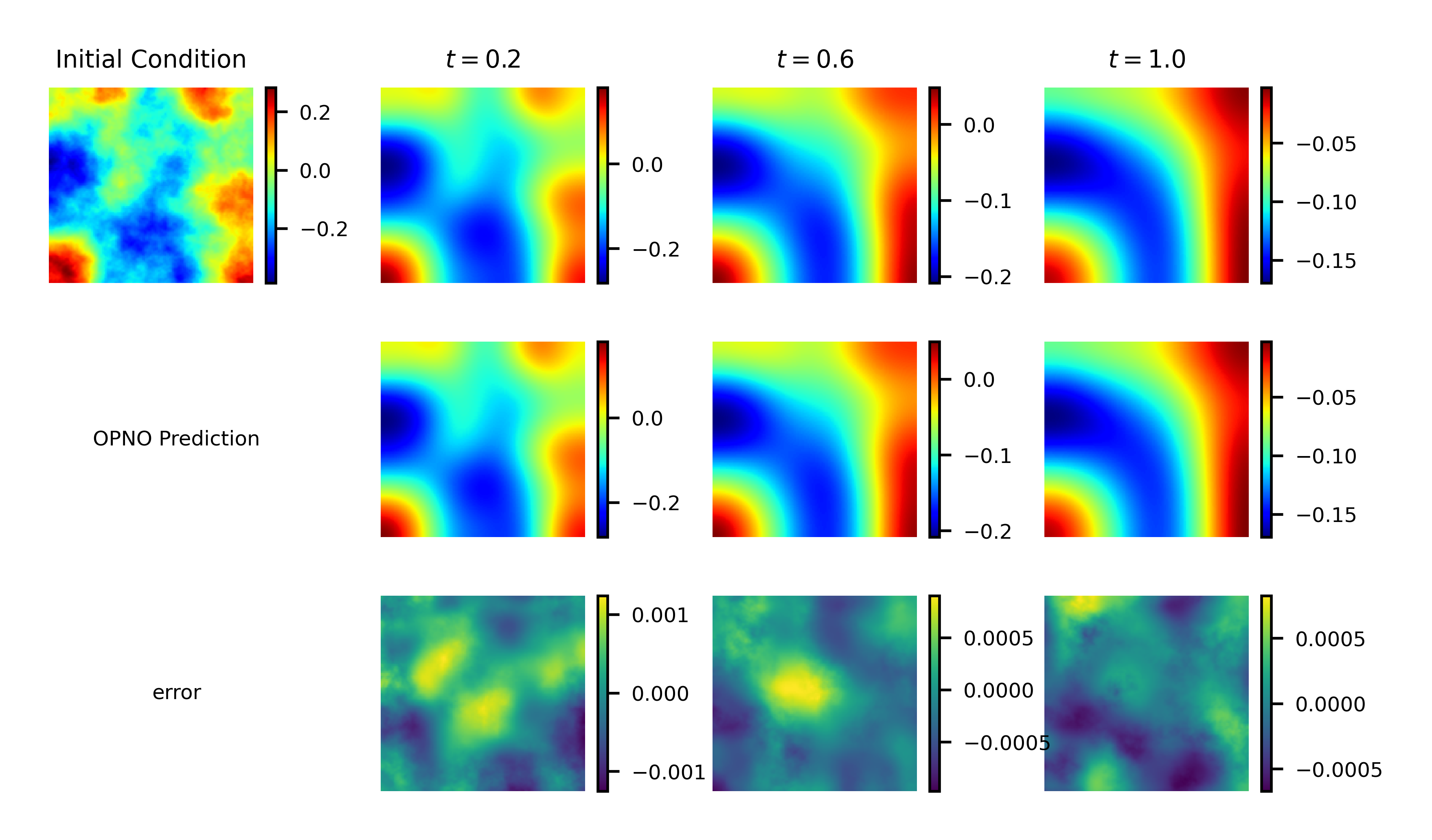}}
  \end{center}
    \caption{\textbf{Two OPNO prediction instances of 2D experiment} that are trained on $51\times 51$ grids and tested on $201\times 201$ grids. For each example, the upper figures are the reference solutions in the test data, while the
middle figures are the predictions of the OPNO at the corresponding time based on the input initial condition, and the bottom row represents the pointwise errors of the OPNO predictions.}
   \label{fig:2d-plots}
\end{figure}

The FNO and OPNO models are trained for $3000$ epochs with an initial learning rate $0.001$ that is halved every $300$ epochs, with their batch size fixed to $20$. In addition, we deploy the neural operators with 16 modes and 24 channel width. \cref{tab:exp-4} shows that the OPNO also outperforms the FNO model in
solving the 2d Burgers' equation, while two instances in test set for prediction is given in  \cref{fig:2d-plots}.
\begin{table}[htbp]
  \footnotesize
  \caption{Benchmark results obtained on the 2D Burgers' equation with Neumann BCs.}
  \label{tab:exp-4}
\begin{center}
  \begin{tabular}{ccc|cc|cc}
    \hline
     N &\multicolumn{2}{c|}{50} &\multicolumn{2}{c|}{100}&\multicolumn{2}{c}{200}\\ \cline{2-7}
  &$L^2$&b.c. $L^{\infty}$ &$L^2$&b.c. $L^{\infty}$ &$L^2$&b.c. $L^{\infty}$ \\\hline
    FNO &$5.278e-3$&$1.625e-1$ & $5.892e-3$ & $3.590e-1$   &$6.913e-3$&$7.849e-1$\\
    \hline
    OPNO &$3.714e-3$&$2.886e-12$ & $3.359e-3$ & $1.955e-12$  & $3.354e-3$&$7.931e-12$\\ \hline
  \end{tabular}
\end{center}
\end{table}
\subsection{Discussions on the generalization errors}
\begin{figure}[htbp]
  \begin{center}
    \subfloat[][loss landscape]{\includegraphics[width=.3\linewidth]{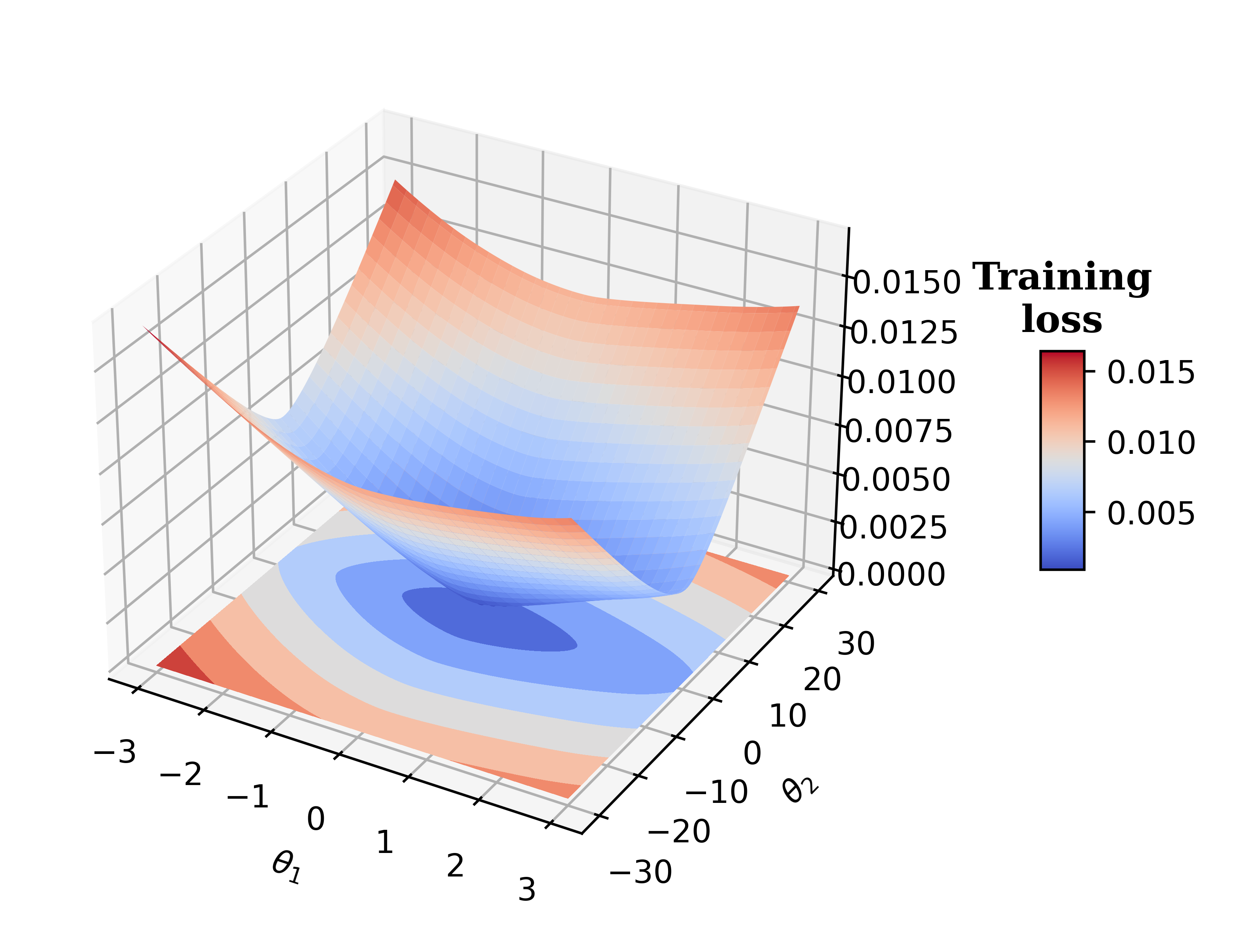}\label{fig:landscape1}}
    \subfloat[][contour plot for training and test error landscapes]{\includegraphics[width=.25\linewidth]{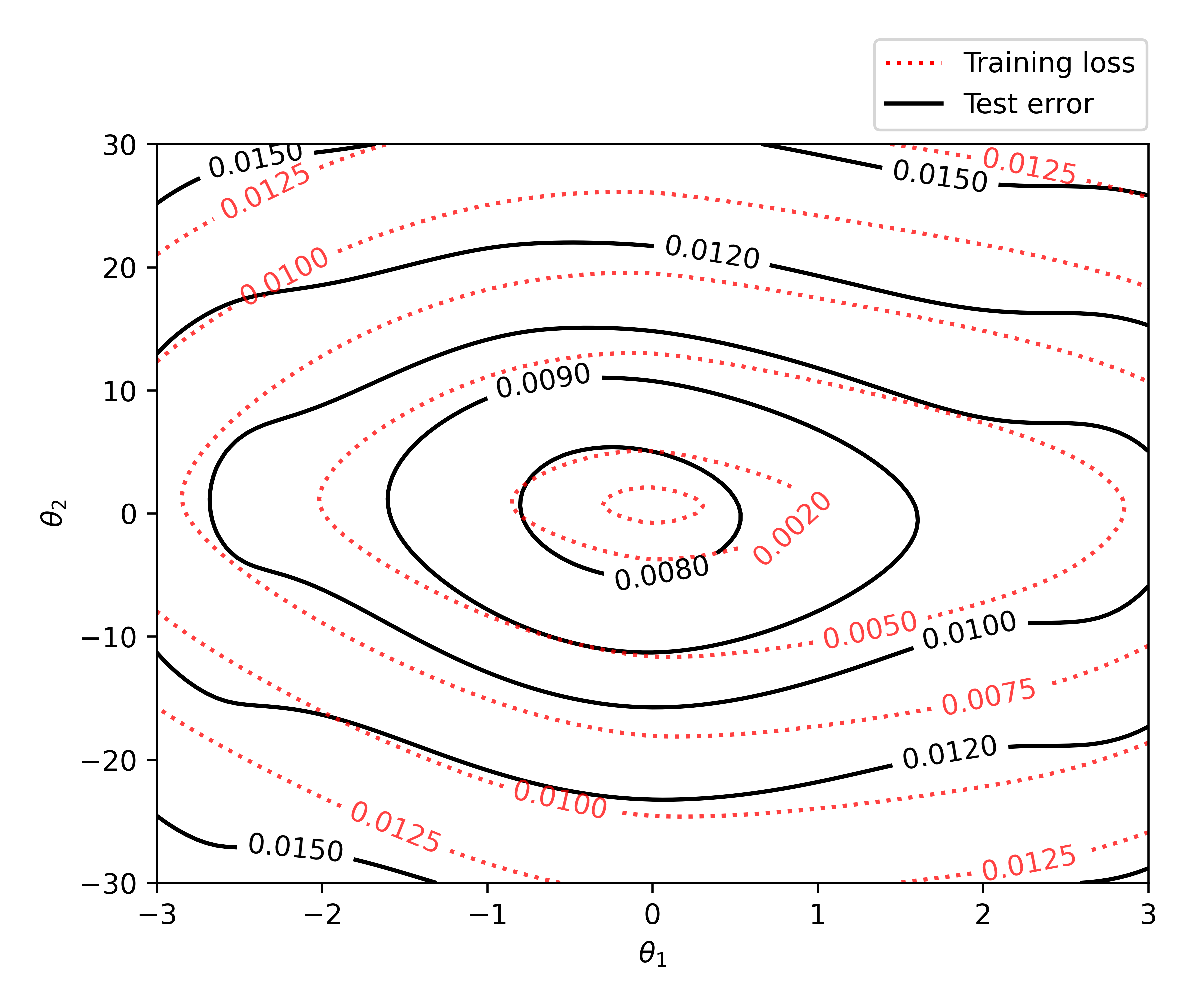}\label{fig:landscape2}}\quad
    \subfloat[][zoomed-in view of the landscapes]{\includegraphics[width=.25\linewidth]{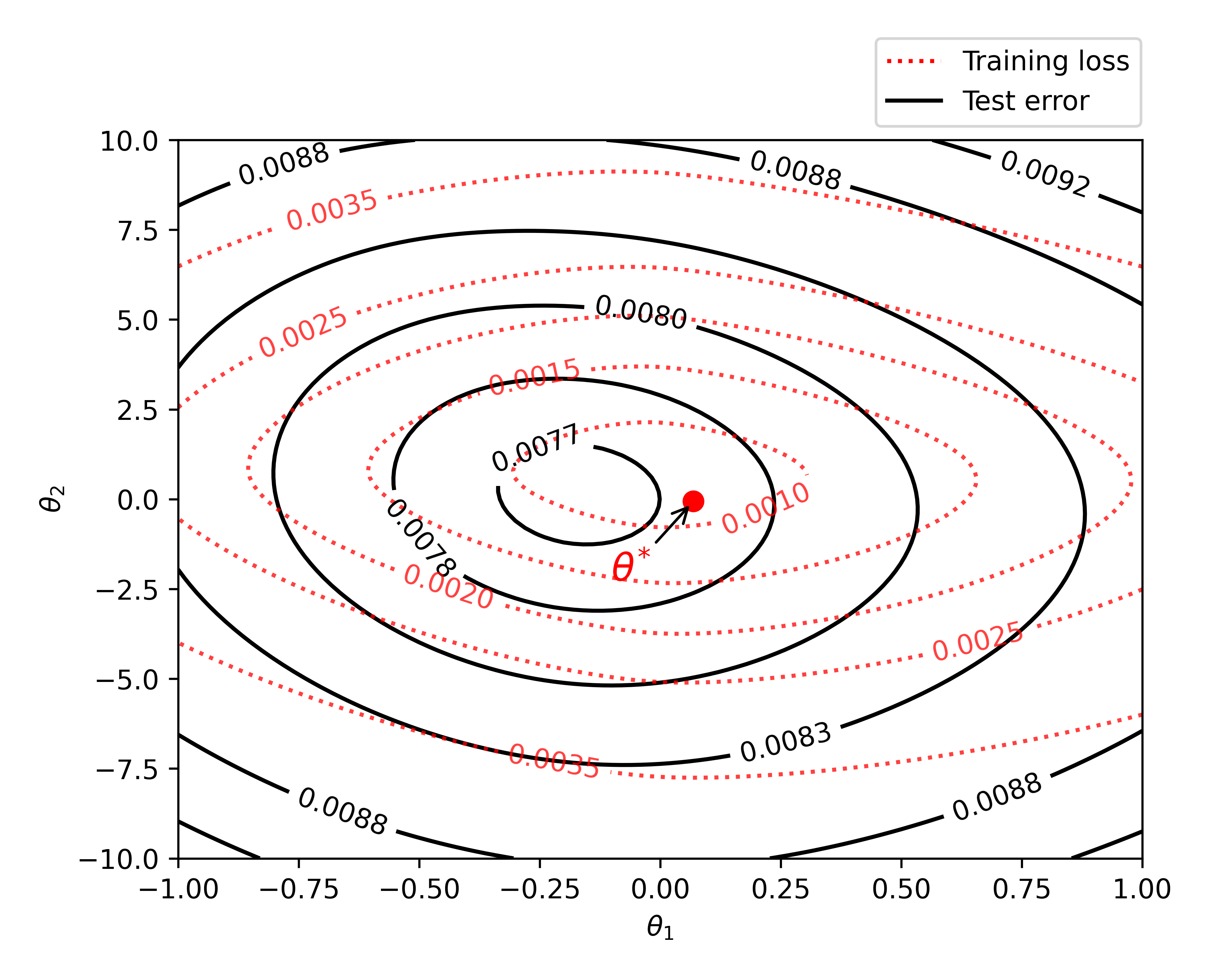}\label{fig:landscape2-zoom}}
  \end{center}
    \caption{The training and test error landscapes of OPNO.}
\end{figure}

Based on the visualization of the training and test error landscapes, we aim to provide an intuitive explanation for the better performance and non-overfitting phenomena of OPNO. The training error landscape of the OPNO model for
the Neumann-BC experiment (\cref{exp:burgers}) is demonstrated in  \cref{fig:landscape1}, where the loss function is almost convex with respect to the learnable parameters $\theta_1$ and $\theta_2$, making the OPNO easy to train. In addition, the
contour plot about the comparison of the training and test error landscapes, along with its zoomed-in subgraph, are demonstrated in \cref{fig:landscape2} and \cref{fig:landscape2-zoom}, respectively. One can see that the geometries of the two landscapes
are pretty close.
  
The similarity between the training and test error landscapes can result in the non-overfitting phenomena. Furthermore, it leads to a supervised model with low theoretical generalization error, which is usually decomposited into two
components \cite{shalev2014understanding}: the ``approximation error'' (the minimum
error achievable by the model in the hypothesis class, introduced by the model bias) and the ``estimation error'' (the difference between the approximation error and the minimum error achievable on
the given data,
depending on the size of the finite training data and the richness of hypothesis class).

While for one model, there is always a trade-off between these two errors (also known as the
bias-variance trade-off), introducing the BC satisfaction architecture, however, can reduce both of them simultaneously. On the one hand, the approximation error is reduced by eliminating the model bias caused by wrong underlying function
basis. On the other hand, the estimation error also decreases because the hypothesis space is limited to the class of operators with output strictly satisfying the BCs. 

Moreover, the  near-convex loss function enables the OPNO to practically obtained the achievable minimum error through long-term training. Consequently, as the result of the Robin-BC experiment (\cref{exp:robin}) shows, if the size of the dataset is sufficiently large, the test errors of BC-satisfying neural
operators can reach an impressively low level in practice.
\section{Conclusion}
In this study, based on the architecture of the Spectral Operator Learning, a boundary-condition-satisfying neural operator named Orthogonal Polynomial Neural Operator is proposed for solving PDEs with
Dirichlet, Neumann and Robin BCs. As far as we are concerned, the OPNO is the first deep-learning method that can
generate solutions strictly satisfying such non-periodic BCs for a family of PDEs, while extensive experiments verify that our model can achieve the state-of-the-art performance in solving PDEs with these BCs. An
analysis of the worst-case performance illuminates the key role of the BC satisfaction property in developing precise
and reliable neural operators. 

Based on a massive dataset, for the first time, we confirm that a deep-learning
method can be more accurate in solving PDEs than a 2nd-order numerical method with a considerably fine mesh while being $4$ to $5$ orders of magnitude faster. The Universal Approximation Theorem is also provided, which ensures the expressive
capacity of the proposed method.

Additionally, the OPNO model possesses other desirable computational properties: It is a fast algorithm with a time complexity of $O(N\log N)$, which is important for big-data training and large-scale
predictions; models that are trained on a coarse mesh can be directly applied on a fine mesh
without loss of numerical accuracy; and the model parameters are sparse and memory-efficient. These properties result from our adpotion of the structure in spectral numerical methods for
SOL architecture.

\appendix
\section{Chebyshev polynomials and their properties}\label{app:chebyshev}
According to the well-known Weierstrass approximation theorem, any continuous function $u(x) \in C(I)$ can be uniformly approximated by a polynomial function. Furthermore, the Chebyshev
polynomials (of the first kind) form one of the most popular polynomials basis  for developing numerical methods. They are given by the three-term recurrence relation
\begin{equation}
T_{n+1}(x) = 2 x T_n(x) - T_{n-1}(x), \ n \geq 1, \ x \in I. \nonumber
\end{equation}
Fortunately, they have an explicit formula
\begin{equation}
 T_n(x) = \cos n \theta, \ \theta = \arccos(x), \ x \in I, \nonumber
\end{equation}
from which we can easily derive the following properties:
\begin{eqnarray}
  &&T_n(\pm 1) = (\pm 1)^n, \ T_n'(\pm 1) = (\pm 1)^{n-1} n^n, \nonumber\\
  &&\int_{-1}^1 T_n(x) T_m(x) \frac{1}{\sqrt{1-x^2}} dx = \frac{c_n \pi}{ 2} \delta_{mn}  \label{eq:ortho-omega}\\ %\ \text{(Orthogonality)},\\
  &&T'(x) =  2n \sum\limits_{k=0 \atop k+n \text{odd}}^{n-1} \frac{1}{c_k} T_k(x),\nonumber
\end{eqnarray}
where $c_0 = 2$ and $c_n = 1$ for $n \geq 1$. Eq. \eqref{eq:ortho-omega} shows that the Chebyshev polynomials are orthogonal with respect to the weight function $\omega(x)
= (1-x^2)^{-\frac{1}{2}}$, and thus for a given $N > 0$, the
$L_{w}^2$-orthogonal projection $\Pi_N: u(x) \mapsto \left\{ \hat u_k \right\}_{0\leq k \leq N}$ is defined as
\begin{equation}
\label{eq:33}
\Pi_N(u) (x) \triangleq \sum\limits_{k=0}^N  \hat u_k T_k(x), \nonumber
\end{equation}
where
\begin{equation}
\hat u_k = \int_{-1}^{1}  u(x) T_k(x) \omega(x) dx. \nonumber
\end{equation}
Remark that the uniqueness of interpolation polynomials yields $I_N u = \Pi_N u, \ \forall u \in C(I)$.

Denote $\mathbb P_{N}(I)$ as the space of all polynomials on $I$ of the degree no greater than $N$. For the CGL quadrature, the Chebyshev-Gauss-type nodes and weights are defined as
$$x_j = - \cos \frac{\pi j}{N}, \ \omega_j = \frac{\pi}{c_j N}, \ 0 \leq j \leq N,$$
where $c_{j}=1$ for $j=1, 2, ..., N-1$ but $c_0 = c_N = 2$. As a consequence, the following equation holds:
\begin{equation}
\int_{-1}^1 p(x) \frac{1}{\sqrt{1-x^2}} dx = \sum\limits_{j=0}^N p(x_j) \omega_j, \ \forall p \in \mathbb P_{2N-1}.   \nonumber
\end{equation}
What's more, the discrete Chebyshev transform
$$\mathcal C_h: \left\{ u(x_j) \right\}_{0 \leq j \leq N} \mapsto \left\{  \hat u_k \right\}_{0 \leq k \leq N}$$
can be carried out with $O(N \log N)$ operations via fast Fourier transform (FFT).

For $u(x) = \sum\limits_{k=0}^N \hat u_k T_k(x) \in \mathbb P_N(I)$, its derivative $u'=\sum\limits_{k=0}^{N-1} \hat u^{(1)}_k T_k(x) $ is obtained by performing differentiation on the frequency space, namely,
\begin{equation}
  \label{eq:poly-diff}
\left\{\begin{aligned}
    &\hat u_{N-1}^{(1)} = 2N \hat u_N, \\
    &\hat u_{k-1}^{(1)} = (2k \hat u_k + \hat u_{k+1}^{(1)}) / \tilde c_{k-1}, \ k = N-1, ..., 1,
\end{aligned}\right.
\end{equation}
which can also be calculated via a linear convolution procedure in $O(N log N)$ operations.

\section{The fast compacting transform for Neumann conditions}\label{app:cp-neumann}
For the compact combination $\phi_k(x) = T_k(x) - \frac{k^2}{(k+2)^2}T_{k+2}$ satisfying the BCs in Eq. \eqref{eq:bc-neumann}, we have
the \textbf{forward compacting transform}
\begin{equation}
\beta_j =\left\{\begin{aligned}
  &\alpha_j, \ j = 0, 1, \\
  &p_{j} \beta_{j-2} + \alpha_{j}, \ 2 \leq j \leq N-4, \\
  &- \frac{1}{p_j} \alpha_{j+2}, \ j = N-3, N-2, \\
\end{aligned}\right. \nonumber
\end{equation}
and the \textbf{backward compacting transform}
\begin{equation}
\alpha_j =\left\{\begin{aligned}
    &\beta_j, \ j = 0, 1, \\
    &\beta_j - p_{j} \beta_{j-2}, \ 2 \leq j \leq N-2, \\
    &- p_{j} \beta_{j-2}, \ j = N-1, N,
\end{aligned}\right. \nonumber
\end{equation}
where $p_{j+2} = j^2 / (j+2)^2, 0 \leq j \leq N-2$.

Similarly, the only troublesome part when conducting parallel computing is the forward transform with $2 \leq j \leq N-2$. However, let $\tilde \alpha_j = j^2\alpha_j, \tilde \beta_j =j^2 \beta_j$; then, we obtain the following equation:
\begin{equation}
\tilde \beta_j = \tilde \beta_j + \tilde a_j , \ 2 \leq j \leq N-2,  \nonumber
\end{equation}
which can be efficiently computed with the linear convolution method. A similar approach can also be applied to the Robin BCs.

\section{Proof of the Universal Approximation Theorem \cref{thm:uat}}\label{sec:app-uat}
The constructive proof of the Universal Approximation Theorem for the FNO has been given in \cite{kovachki2021universal} by Nikola Kovachki, Samuel Lanthaler, and Siddhartha Mishra. %The proof of
                                %lemma \cref{lem:conv1} is almost identical to that in \cite{kovachki2021universal}, the major change being the substitution of the Shen transform for the Fourier
                                %transform. We will be dedicated to the concrete construction of the desired OPNO.
Besides, Hao Liu, Haizhao Yang, et al. also estimated the generalization error of neural operators with encoders and decoders of trigonometric functions or Legendre polynomials in \cite{haizhao2022legendre}. Actually, Jacobi's expansions \cite{mason2002chebyshev} for the underlying basis of FNO show that
\begin{eqnarray}
\label{eq:35}
  &&\sin(k\pi x) = 2 \sum\limits_{n=0}^{\infty}(-1)^n J_{2n+1}(k\pi)T_{2n+1}(x),  \nonumber\\
  &&\cos(k\pi x) = J_0(z) + 2 \sum\limits_{n=1}^{\infty}(-1)^n J_{2n}(k\pi) T_{2n}(x),  \nonumber
\end{eqnarray}
where $J_n(x)$ is the Bessel function of the first kind. Using the asymptotic formula
\begin{equation}
\label{eq:48}
J_{k}(x) \sim \frac{1}{2k\pi} \bigg( \frac{ex}{2k} \bigg)^k, \ k \gg 1, \nonumber
\end{equation}
we find that the expansion coefficients exponentially decay as long as $N$ is sufficiently large. Therefore, one can guess from intuition that it should be feasible to extend the Universal Approximation
Theorem to the OPNO model with finite modes and bandwidth.

Denote the $L^2$ norm $\left\| \cdot \right\|_{L^2}$ by $\left\| \cdot
 \right\|$ in short. Before we start the proof, as a result of the following lemma given in \cite{kovachki2021universal}, we adjust the direction of the proof to approximating the ``spectral projection'' operator $\mathcal G_N \triangleq \Pi_N \mathcal G (\Pi_N a)$ instead of the original operator $\mathcal G$
using neural operators:

\begin{lemma}[\cite{kovachki2021universal}]
  \label{lem:conv1}
$\forall \epsilon > 0$, there exists $N \in \mathbb N$ such that

\begin{equation}
\label{eq:GN}
  \left\| \mathcal{G} (a) - \mathcal{G}_N(a) \right\| \leq \epsilon, \forall a \in K. \nonumber
\end{equation}
\end{lemma}
\cref{lem:conv1} is a consequence of the following facts: Since $K \subset H^s(I^d)$ is compact and $\mathcal G$ continuous, the set $\hat K \triangleq K \cup \bigcup\limits_{N \in \mathbb{N}} \Pi_N K$ and its image
$\mathcal G(\hat K)$ are also compact, and $\mathcal G \big|_{\hat K}$ uniformly continuous. The compactness of $\mathcal G(\hat K)$ gives the existence of a modulus of continuity $\omega : [0,
\infty) \rightarrow [0, \infty)$, such that
\begin{equation}
  \left\| \mathcal G(a) - \mathcal G(a')  \right\| \leq \omega(\left\| a - a' \right\|_{H^s} ), \ \forall a, a' \in \hat K, \nonumber
\end{equation}
while the compactness of both $\hat K$ and $\mathcal{G}(\hat K)$ yeilds
\begin{equation}
\limsup\limits_{N \to \infty} \sup\limits_{v \in \mathcal G(\hat K)} \left\| (1 - \Pi_N) v \right\| = 0 = \limsup\limits_{N \to \infty} \sup\limits_{a \in \hat K} \left\| (1-\Pi_N) a \right\|_{H^s}.
\end{equation}
So we have the inequality 
\begin{equation}
\begin{split}
  \left\| \mathcal{G} (a) - \mathcal{G}_N(a) \right\| &= \left\|   \mathcal G(a) - \Pi_N \mathcal G(\Pi_N a)  \right\| \\
                                                            &\leq \left\| \mathcal G(a) - \Pi_N \mathcal G(a) \right\| + \left\| \Pi_N \mathcal G(a) - \Pi_N \mathcal G(\Pi_N a) \right\|\\
                                                            &\leq \sup\limits_{v \in \mathcal G(\hat K)} \left\| (1-\Pi_N) v \right\| + \omega(\sup\limits_{a \in \hat K} \left\| a-\Pi_N a
                                                               \right\|_{H^s} )\\
                                                              & \leq \epsilon.
\end{split}
\nonumber
\end{equation}
holds once $N$ is sufficiently large.

% Then, a simple argument shows that both N-indexed sequences on the RHS of the following triangle inequality vanish at infinity:
 
Thus, the problem reduces to constructing an OPNO that approximates $\mathcal G_N$. After defining an auxiliary operator $\hat {\mathcal G}_N $ such that
\begin{equation}
\label{eq:36}
\mathcal G_N(a) = \mathcal{C}_h^{-1} \circ \hat {\mathcal G}_N \circ [\mathcal{C}_h \circ \Pi_N](a), \nonumber
\end{equation}
we now perform its individual steps, in which we set $d_a = d_v = d = 1$ without any loss of generalization.

\subsection{Approximation of $\mathcal C_h \circ \Pi_N$ (the discrete Chebyshev transform)}
\begin{figure}[tbhp]
  \centerline{\includegraphics[width=0.75\textwidth]{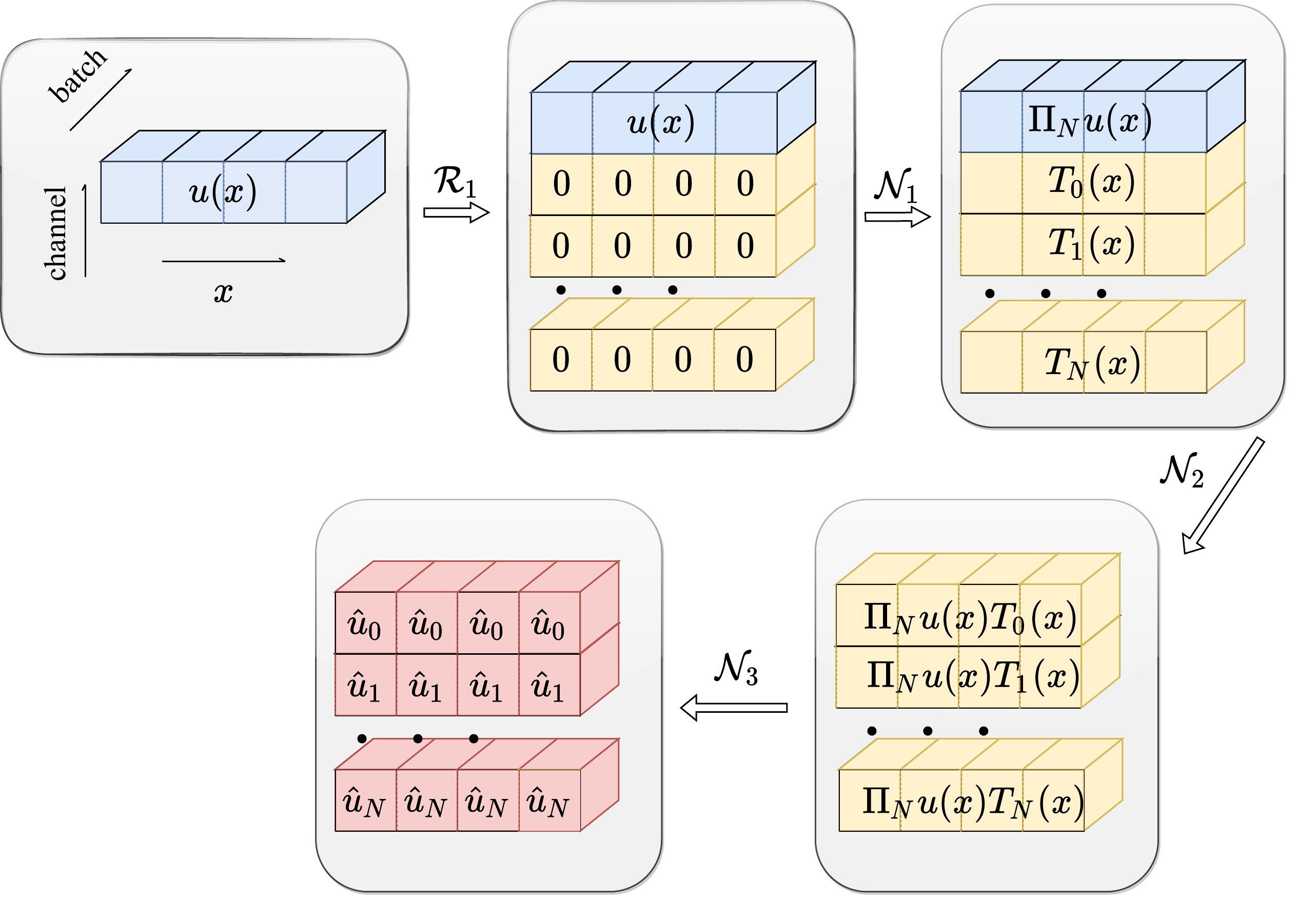}}
  \caption{Schematic diagram for the OPNO approximation of the operator $\mathcal C_h \circ \Pi_N$}
  \label{fig:cubes-for-uat}
\end{figure}

The trick used in the proof is to manually build a discrete Chebyshev polynomial decomposition (or its inverse in the next subsection) in the \textbf{channel} dimension, even though in the structure
of the OPNO, such a transform is performed in the spatial dimension. Please refer to  \cref{fig:cubes-for-uat} for the schematic diagram.
Specifically, instead of an ($N+1$)-dimensional vector, we view the output of the target operator as a \textbf{constant ($N+1$)-dimensional vector-valued function}.
In other word, we aim to approximate the operator
\begin{equation}
\label{eq:46}
\mathcal C_h \circ \Pi_N : L^2(I; \mathbb R) \to L^2(I; \mathbb R^{N+1}), \ v(x) \mapsto \hat w(x) \equiv \left\{ \hat v_k\right\}_{k \in K_N}, \nonumber
\end{equation}
where $K_N = \left\{ {0, 1, 2, ..., N} \right\}$ and $\hat v_k$ is the $k$-th Chebyshev polynomial expansion coefficient of $v(x)$.

It is notable that the OPNO layer (Eq. \cref{eq:25}) reduces to a vanilla (ResNet) neural network when filling all parameters of $A_l$ with $0$. In other words, feedforward neural networks
can be seen as trivial instances of the OPNO. So we denote by $\mathbf a \vee \mathbf b$ the
concatenation in the channel dimension of the two vectors herein, and define a lifting
operator $\mathcal R_1$ as well as two continuous operators $\mathcal N_1, \tilde {\mathcal N}_2$ such that
\begin{eqnarray}
\mathcal R_1 &:& L^2(I) \to L^2(I; \mathbb R^{N+2}), \ v(x) \mapsto w^{(0)}(x) \triangleq \left\{ v(x) \right\} \vee \mathbf{0}_{N+1}, \nonumber\\
\mathcal N_1 &:& L^2(I; \mathbb R^{N+2}) \to L^2(I; \mathbb R^{N+2}), \ w^{(0)}(x) \mapsto w^{(1)}(x) \triangleq \left\{ \Pi_N w^{(0)}_{0}\right\} \vee \left\{ T_k(x) \right\}_{k \in K_N},\nonumber\\
\tilde{\mathcal N}_2 &:& L^2(I; \mathbb R^{N+2}) \to L^2(I; \mathbb R^{N+1}), \ w^{(1)}(x) \mapsto w^{(2)}(x) \triangleq \left\{ w^{(1)}_0 \cdot w^{(1)}_{1+k} \right\}_{k \in K_N}.\nonumber
\end{eqnarray}
Therefore,
\begin{equation}
\label{eq:41}
\tilde{\mathcal N}_2 \circ \mathcal N_1 \circ \mathcal R_1(v) (x) = \left\{ \Pi_N v(x) \cdot T_k(x) \right\}_{k \in K_N}. \nonumber
\end{equation}
Notably, $\mathcal N_1(w^{(0)})(x) = b_1(x) + \mathcal S^{-1} \big( \mathcal A_1 \mathcal C_h w^{(0)} \big) (x)$ is an OPNO layer with

$$\ b_1(x) = 0 \vee \left\{ T_k(x) \right\}_{k \in K_N},  \ \mathcal A_1 =  \mathcal C_h \circ \mathbf 1_{k \leq N}, $$
and $\tilde{\mathcal N}_2$ is simply a multiplication operation that can be easily approximated by a $\sigma$-activated neural network $\mathcal N_2$.

Furthermore, using the fact that
\begin{equation}
\label{eq:42}
\mathcal C_h[w^{(2)}(\cdot)](j) =  \mathcal C_h \bigg[\Pi_N v(\cdot) T_k(\cdot)\bigg] (j) = \int_{-1}^1 \Pi_N v(x) T_k(x) T_j(x) w(x) dx \nonumber
\end{equation}
and plugging $j=0$ back into the above equation, we find that
\begin{equation}
\label{eq:43}
\mathcal C_h \bigg[\Pi_N v(\cdot) T_k(\cdot)\bigg] (0) = \la \Pi_N v(x), T_k(x) \ra_w = \hat v_k. \nonumber
\end{equation}
In other words, letting $\mathcal A_3 = \mathcal C_p \circ \mathbf 1_{k=0}$, we have
\begin{equation}
\label{eq:34}
\mathcal S^{-1} \mathcal A_3 \mathcal C_h [\Pi_N v(\cdot) T_k(\cdot)] = \mathcal C_h^{-1}  \left\{ \hat v_k \vee \mathbf 0_{N} \right\} = \hat v_k T_0(x) \equiv \hat v_k. \nonumber
\end{equation}
Consequently, we define an OPNO layer $\mathcal N_3$ such that
\begin{equation}
\label{eq:45}
\mathcal N_3 : L^2(I; \mathbb R^{N+1}) \to L^2(I; \mathbb R^{N+1}), \  \tilde w^{(2)}(x) \mapsto w^{(3)}(x) \triangleq \mathcal S^{-1} \mathcal A_3 \mathcal C_h w^{(2)}. \nonumber
\end{equation}

Finally, let $\mathcal L_1 = \mathcal N_3 \circ \mathcal N_2 \circ \mathcal N_1 \circ \mathcal R_1$, then $\mathcal L_1$ is the desired answer. Broadly speaking, $\mathcal L_1$ is an OPNO ``up to $\sigma$-activation'' since the activation function is missing in $\mathcal N_3$. Such a difference, however, is insignificant due to the flexibility of neural
networks, especially in practical applications, and can always be fixed by composition with a $\sigma$-activated neural network approximating the identity.

\subsection{Approximation of $\mathcal C_h^{-1}$}
The operator $\mathcal C_h^{-1}$ is seen as an operator with a constant vector-valued function $\hat w(x) \in  L^2(I; \mathbb R^{N+1})$ as input, where
$$\hat w(x)\equiv \left\{ \hat w_k \right\}_{k \in K_N}.$$
Let
\begin{eqnarray}
\label{eq:37}
&&\mathcal R_{2} : L^2(I; \mathbb R^{N+1}) \to L^2(I; \mathbb R^{2N+2}), \ \hat w(x) \mapsto \hat w^{(0)} \triangleq \hat w(x) \vee \left\{T_k(x) \right\}_{k \in K_N},  \nonumber\\
&& \tilde {\mathcal N}_4 : L^2(I; \mathbb R^{2N+2}) \to L^2(I; \mathbb R), \ \hat w^{(0)}(x) \mapsto w(x) \triangleq \sum_{k=0}^N \hat w^{(0)}_k(x) \cdot w^{(0)}_{k+N+1}(x). \nonumber
\end{eqnarray}
Then, we have
\begin{equation}
\label{eq:53}
\tilde {\mathcal N}_4 \circ \mathcal{R}_2 (\hat w(x)) = \sum\limits_{k=0}^N  \hat w_k T_k(x) = \mathcal C_h^{-1} \hat w(x). \nonumber
\end{equation}
Apparently, $\mathcal R_2(\hat w(x)) = W_4 w(x) + b(x)$ with
$$W_4 = \begin{bmatrix}
\mathbf E_{(N+1) \times (N+1)} \\ \mathbf 0_{(N+1) \times (N+1)}
\end{bmatrix} \in \mathbb R^{(2N+2) \times (N+1)}, \ b(x) = \mathbf 0_{N+1} \vee \left\{ T_k(x) \right\}_{k \in K_N},$$
and there exists a conventional neural network $\mathcal N_4$ approximating $\tilde {\mathcal N}_4$ to any desired
accuracy. Consequently, the OPNO $\mathcal L_2 \triangleq \mathcal N_4 \circ \mathcal R_{2}$ is the desired answer.

\subsection{Approximation of $\hat {\mathcal{G}}_N$}\label{app:uat-3}
It remains to be shown that the operator $\hat {\mathcal{G}}_N$ can be approximated by an OPNO, which is straightforward since this operator is nothing but a continuous map between
compact subsets of
\textbf{finite}-dimensional function spaces. According to the Universal Approximation Theorem of feedforward neural networks, there exists a (degenerate) OPNO $\mathcal L_3$ that approximates
$\mathcal G_N$ to an arbitrary level of precision. Select the activation function $\sigma$ as a Lipschitz function, then $\mathcal L_3$ is thus also Lipschitz continuous with a Lipschitz constant $\text{Lip}(\mathcal
 L_3)$.
 
Remind that the image of a compact set remains compact under any continuous function (operator), and any continuous function on a compact set attains its maximum. Following the processes above, for any fixed $\epsilon >0$, we
 construct the OPNOs $\mathcal L_1, \mathcal L_2, \mathcal L_3$ such that the following equations hold for any compact set $K \subset H^s(I)$:
\begin{equation}
\begin{split}
  \text{Lip}(\mathcal L_3) \sup\limits_{a \in K} &\left\| (\mathcal C_h \Pi_N - \mathcal L_1)a \right\| < \frac{\epsilon}{3\left\| \mathcal C_p^{-1} \right\|}  , \\
  \sup\limits_{w \in \mathcal C_h \Pi_N(K)} &\left\| (\hat{\mathcal G}_N  - \mathcal L_3)w \right\| < \frac{\epsilon}{3\left\| \mathcal C_p^{-1} \right\|} , \\
    \sup\limits_{v \in \mathcal L_3 \circ \mathcal L_1(K)} &\left\| (\mathcal C_h^{-1} - \mathcal L_2) v \right\|  < \frac{\epsilon}{3},
\end{split}
\nonumber
\end{equation}
where
$$\left\| \mathcal C_p^{-1} \right\| = \sup\limits_{\hat w(x) \equiv \left\{ \hat w_k \right\}_{k \in K_N}} \left(  \left\| \mathcal C_p^{-1} \hat w(x)\right\| / \left\| \hat w(x) \right\| \right).$$

To summarize, let $\mathcal L = \mathcal L_2 \circ \mathcal L_3 \circ \mathcal L_1$, then
%then since $\epsilon>0$ is arbitrary, we have
\begin{equation}
\begin{split}
  \sup\limits_{a \in K} \left\| \mathcal G_N a - \mathcal L a\right\| &=  \sup\limits_{a \in K}\left\| \mathcal C_h^{-1}\circ \hat{\mathcal G}_N \circ \mathcal C_h\Pi_N a - \mathcal L_2 \circ \mathcal L_3 \circ
                                                                                \mathcal L_1 a\right\| \\
  &\leq \sup\limits_{a \in K} \left\| \mathcal C_h^{-1}\circ \hat{\mathcal G}_N \circ \mathcal C_h\Pi_N a -  \mathcal C_h^{-1} \circ\mathcal L_3 \circ \mathcal L_1 a\right\|\\ 
    & \quad + \sup\limits_{a \in K} \left\| \mathcal C_h^{-1} \circ\mathcal L_3 \circ \mathcal L_1 a - \mathcal L_2 \circ \mathcal L_3 \circ \mathcal L_1a \right\|\\
    & \leq \left\| \mathcal C_p^{-1} \right\|\sup\limits_{a \in K} \left\| \hat{\mathcal G}_N \circ \mathcal C_h \Pi_N a - \mathcal L_3 \circ \mathcal L_1 a \right\| \\
  &\quad + \sup\limits_{v \in \mathcal L_3 \circ \mathcal L_1(K)}
    \left\| \mathcal C_h^{-1} v - \mathcal L_2 v\right\| \\
   &\leq \left\| \mathcal C_p^{-1} \right\| \bigg( \sup\limits_{a \in K} \left\| \hat{\mathcal G}_N \circ \mathcal C_h \Pi_Na - \mathcal L_3 \circ \mathcal C_h \Pi_N a \right\| \\
& \quad + \sup\limits_{a \in K}\left\| \mathcal L_3 \circ \mathcal C_h
     \Pi_N a - \mathcal L_3 \circ \mathcal L_1 a \right\| \bigg) \\  
 & \quad + \sup\limits_{v \in \mathcal L_3 \circ \mathcal L_1(K)} \left\| (\mathcal C_h^{-1} - \mathcal L_2) v \right\| \\
  & \leq \left\| \mathcal C_p^{-1} \right\| \sup\limits_{w \in \mathcal C_h \Pi_N(K)} \left\| (\hat{\mathcal G}_N  - \mathcal L_3)w \right\| \\
  & \quad +  \left\| \mathcal C_p^{-1} \right\| \text{Lip}(\mathcal L_3) \sup\limits_{a \in K} \left\| (\mathcal C_h \Pi_N - \mathcal L_1)a
    \right\| \\
    & \quad + \sup\limits_{v \in \mathcal L_3 \circ \mathcal L_1(K)} \left\| (\mathcal C_h^{-1} - \mathcal L_2) v \right\| \\
    & < \epsilon,
\end{split}
\nonumber
\end{equation}
which complete the proof.

\bibliographystyle{siamplain}
\bibliography{references}

\end{document}